# CHARACTERIZING MARKOV EQUIVALENCE CLASSES FOR AMP CHAIN GRAPH MODELS[1]

BY STEEN A. ANDERSSON AND MICHAEL D. PERLMAN

*Indiana University and University of Washington*

Chain graphs (CG) (= adicyclic graphs) use undirected and directed edges to represent both structural and associative dependences. Like acyclic directed graphs (ADGs), the CG associated with a statistical Markov model may not be unique, so CGs fall into Markov equivalence classes, which may be superexponentially large, leading to unidentifiability and computational inefficiency in model search and selection. It is shown here that, under the Andersson–Madigan–Perlman (AMP) interpretation of a CG, each Markov-equivalence class can be uniquely represented by a single distinguished CG, the AMP essential graph, that is itself simultaneously Markov equivalent to all CGs in the AMP Markov equivalence class. A complete characterization of AMP essential graphs is obtained. Like the essential graph previously introduced for ADGs, the AMP essential graph will play a fundamental role for inference and model search and selection for AMP CG models.

**1. Introduction.** In a graphical Markov model, the nodes of the graph represent the variables of a multivariate statistical distribution, while the edges represent possible dependences. Chain graphs (CG), which may have both undirected and directed edges but no semi-directed cycles, were introduced by Lauritzen and Wermuth [19] and Frydenberg [15] to represent dependences that may be both associative and directional. Cox [12] stated that chain graphs represent "a minimal level of complexity needed to model empirical data." Also see [2, 8, 9, 13, 15, 17, 18, 19, 25, 26, 28, 29].

The LWF Markov property for CGs is an extension of the Markov properties of both acyclic directed graphs (ADG ≡ DAG) and undirected graphs (UG). Recently, Andersson, Madigan and Perlman [1, 4] proposed an *alternative Markov property* (AMP) for CGs that also extends the ADG and UG

Received June 2004; revised June 2005.
[1]Supported in part by NSF Grants DMS-00-71818 and DMS-00-71920.

properties, but that more closely retains the recursive character of ADG models; see [14, 20]. Furthermore, AMP Markov equivalence of CGs (see below), as for ADGs, is determined by their *triplexes*, which have three vertices, while LWF Markov equivalence of CGs is determined by their *complexes*, which have arbitrarily many vertices.

Like ADGs, different CGs may be *Markov equivalent*, that is, may represent the same set of conditional independences (CI), hence, the same statistical models. Because Markov equivalence classes can be superexponentially large even for ADGs (cf. [3]), for the sake of computational efficiency, CG model search ideally should be carried out in the space of CG Markov equivalence classes rather than the space of all CGs.

For any CG $G$, Frydenberg [15] showed the existence of a unique *largest* (i.e., having the most undirected edges) CG $G_\infty$ in the LWF Markov equivalence class containing $G$. Studený [23, 24] proposed that $G_\infty$ be used as a unique representative for the LWF equivalence class. Characterizations of $G_\infty$ have been obtained by Volf and Studený [27] and Roverato [22].

The *ADG essential graph* $D^*$ that uniquely represents the ADG Markov equivalence class of an ADG $D$ was introduced by Andersson, Madigan and Perlman [3]: $D^*$ has the same skeleton as $D$, and contains an arrow $a \to b$ iff this arrow occurs in every member of the equivalence class, whereas it contains a line $a - b$ iff $a \to b$ and $a \leftarrow b$ occur in two different ADGs in the equivalence class. Applications of the ADG essential graph for ADG model search are presented in [10] and [21].

This paper is devoted to the characterization of a unique representative of the AMP Markov equivalence class (temporarily denoted by $[G]$) for a general CG $G$. Andersson, Madigan and Perlman [4] suggested the following extension of the definition of the ADG essential graph: the *AMP essential graph*, temporarily denoted by $G^*$, has the same skeleton as $G$ and contains an arrow $a \to b$ iff this arrow occurs with the same orientation in at least one $G' \in [G]$, but with the opposite orientation in no $G'' \in [G]$. The arrows in $G^*$ are called *essential arrows*. For example, if $G = a \to b - c$, then $[G]$ consists of $G$, $G' := a - b \leftarrow c$, and $G'' := a \to b \leftarrow c$, so $G^* = G''$.

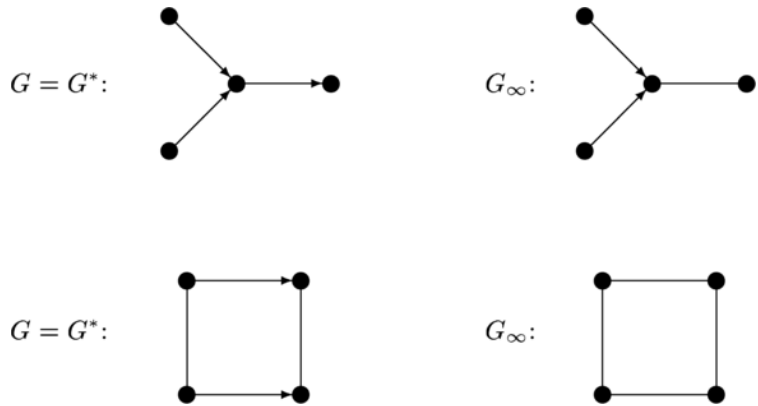

FIG. 1. *Two CGs $G$ for which $G = G^* \neq G_\infty$.*

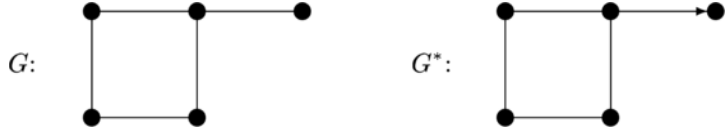

FIG. 2. *An undirected CG $G$ whose essential graph $G^*$ contains a directed edge.*

Because $G_\infty$ is defined for LWF Markov equivalence rather than AMP equivalence, it need not agree with $G^*$, even when $G = D$, an ADG. For example, if $G$ is the first CG (an ADG) in Figure 1, then $G = G^*$, but $G_\infty$ replaces one arrow $\to$ by a line $-$. Similarly, for the second CG $G$ in Figure 1, two arrows are replaced by lines in $G_\infty$. In fact, even if $G$ (and therefore $G_\infty$) is a fully undirected graph, $G^*$ may possess essential arrows—see Figure 2.

In Section 3 we show that the AMP essential graph $G^*$ does in fact uniquely represent its AMP Markov equivalence class $[G]$: $G^*$ is itself a CG (adicyclic) and $G^* \in [G]$ (Theorem 3.2). Section 4 completely describes the configurations of arrows in $G^*$. It is also shown there that if $G$ is itself an ADG $D$ or is AMP Markov equivalent to $D$, then $G^* = D^*$ (Proposition 4.2), and a characterization of those directed graphs that can occur as AMP essential graphs is given (Theorem 4.1). In Section 5 a complete characterization is obtained for AMP essential graphs (Theorem 5.1). Additional results on the structure of AMP essential graphs may be found in [5, 6]. Current research, including an algorithm for constructing $G^*$ from $G$, is reviewed in Section 6.

**2. Graphical terminology.** We write $G \equiv (V, E)$ to indicate a graph $G$ with vertex set $V$ and edge set $E$. Definitions of the graphical terminology and notation used here can be found in [3, 4] and especially [5]. The terms *parent*, *neighbor*, *immorality*, *flag*, *triplex* and *biflag* are particularly important for our study of chain graphs.

Let $a, b$ be distinct vertices of $G$. We write $a \Rightarrow b \in G$ ($a \leftrightarrow b \in G$) ($a \cdots b \in G$) to indicate that either $a \to b \in G$ or $a - b \in G$ (either $a \leftarrow b \in G$ or $a \to b \in G$) (either $a \leftarrow b \in G$, $a - b \in G$, or $a \to b \in G$). A *path* $\pi$ of length $k \geq 1$ from $a$ to $b$ in $G$ is a sequence of distinct vertices $(a \equiv v_0, v_1, \ldots, v_k \equiv b)$ such that $v_{i-1} \Rightarrow v_i \in G$ for all $i = 1, \ldots, k$. A *k-cycle* $(v_0, \ldots, v_k)$ (or simply *cycle*) in $G$ is a path of length $k \geq 3$ with $v_0 = v_k$. A path or cycle $(v_0, \ldots, v_k)$ is *undirected* if $v_{i-1} - v_i \in G$ for all $i = 1, \ldots, k$; it is *directed* (*semi-directed*) if $v_{i-1} \to v_i \in G$ for all (at least one) $i = 1, \ldots, k$. A *chain graph* (CG) is an *adicyclic* graph, that is, has no semi-directed cycles. An induced subgraph of an adicyclic graph is adicyclic. The set of *chain components* of the CG $G$, denoted by $\Xi \equiv \Xi(G)$, is the set of connected components obtained by removing all arrows from $G$.

A path (cycle) $(v_0, \ldots, v_k)$ in $G$ with $k \geq 2$ ($k \geq 4$) is *chordless* if no two nonconsecutive [nonconsecutive (mod $k$)] vertices are adjacent. A *chordless 2-dipath* in $G$ is an induced subgraph of the form $a \to b \to c$. An *antiflag* is an induced subgraph of the form $a - b \to c$.

An undirected graph $G$ is *chordal* ($\equiv$ *decomposable*) if it contains no chordless cycles. The edges of a chordal UG $G \equiv (V, E)$ can be assigned a *perfect orientation* (i.e., acyclic with no immoralities) by the maximum cardinality search (MCS) algorithm (cf. [7], Theorem 2.5, [11], Chapter 4.4). MCS begins by assigning the number 1 to an arbitrary vertex of $G$, then assigning the numbers $2, \ldots, |V|$ consecutively to the remaining vertices, each time selecting the vertex with the most previously numbered neighbors in $G$, breaking ties arbitrarily. The edges of $G$ are then oriented in accordance with this numbering. This numbering is called *perfect* because this orientation can be shown to be perfect. Furthermore, if $A \subseteq V$ is complete, MCS can begin at any $v \in A$ and visit all vertices in $A$ before visiting any vertex in $V \setminus A$. Thus, any edge $a - v \in G$ with $a \in A$ and $v \in V \setminus A$ can be oriented by MCS as $a \to v$.

**3. The essential graph for an AMP chain graph model.** In this paper $G_0 \equiv (V, E_0)$ shall denote a fixed but arbitrary chain graph and $\mathcal{G}$ its AMP Markov equivalence class, that is, the collection of all CGs $G \equiv (V, E)$ such that $\mathcal{P}(G) = \mathcal{P}(G_0)$, where $\mathcal{P}(G)$ is the set of all multivariate probability distributions that satisfy the AMP Markov property specified by $G$. AMP Markov equivalence was characterized by Andersson, Madigan and Perlman [4] as follows.

THEOREM 3.1. *Two chain graphs with the same vertex set are Markov equivalent iff they have the same skeleton and the same triplexes.*

We shall show that $\mathcal{G}$ is uniquely represented by its *essential graph* $\mathcal{G}^*$, defined below. To emphasize that the AMP essential graph depends on $\mathcal{G}$, we now denote it by $\mathcal{G}^*$, rather than by $G_0^*$ as in Section 1.

DEFINITION 3.1. The AMP *essential graph* $\mathcal{G}^* \equiv (V, E^*)$ determined by $\mathcal{G}$ is a graph with the same skeleton as $G_0$. An arrow $a \to b$ occurs in $\mathcal{G}^*$ iff $a \to b$ occurs in at least one $G \in \mathcal{G}$, but $a \leftarrow b$ occurs in no $G' \in \mathcal{G}$. A line $a - b$ occurs in $\mathcal{G}^*$ iff either $a - b \in G$ for every $G \in \mathcal{G}$ or there exist $G, G' \in \mathcal{G}$ such that $a \to b \in G$ and $a \leftarrow b \in G'$.

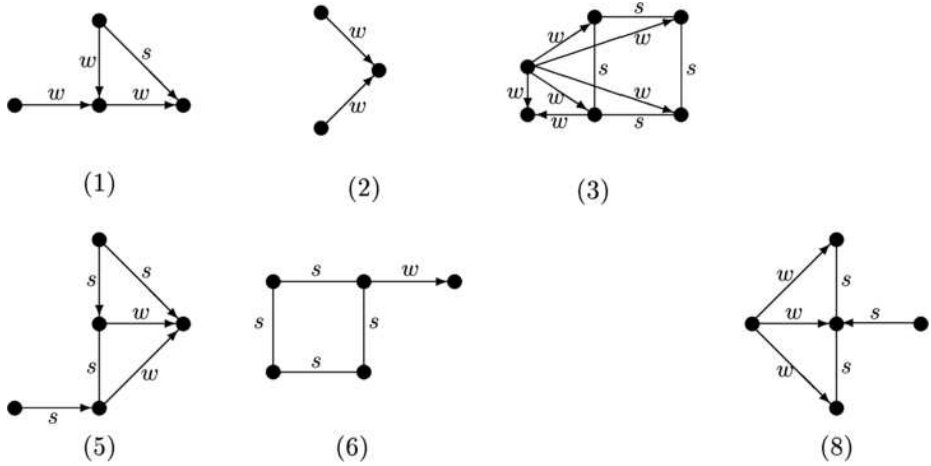

FIG. 3. *Strong/weak arrows/lines in AMP essential graphs (cf. Figure 6).*

Thus, the line $a - b \in \mathcal{G}^*$ iff *either*: $a - b \in G$ for all $G \in \mathcal{G}$, in which case it is called a *strong line* and denoted as $a \overset{s}{-} b$, *or*: $a \to b \in G$ and $a \leftarrow b \in G'$ for some $G, G' \in \mathcal{G}$, in which case it is called a *weak line* and denoted as $a \overset{w}{-} b$. In these two cases, $a$ is a *strong* (resp., *weak*) *neighbor* of $b$. An arrow $a \to b \in \mathcal{G}^*$ is called *strong* and denoted as $a \overset{s}{\to} b$ if it occurs in each $G \in \mathcal{G}$; otherwise it is called *weak* and denoted as $a \overset{w}{\to} b$, in which case $\exists\, G, G' \in \mathcal{G}$ such that $a \to b \in G$ and $a - b \in G'$, while $a \Rightarrow b \in G''$ for all other $G'' \in \mathcal{G}$. In these two cases, $a$ is a *strong* (resp., *weak*) *parent* of $b$. The set of strong (weak) parents in $\mathcal{G}^*$ of a subset $B \subset V$ is denoted by $\mathrm{sp}_{\mathcal{G}^*}(B)$ [$\mathrm{wp}_{\mathcal{G}^*}(B)$].

In the AMP essential graph $a \to b \leftarrow c$, both arrows are weak, while in $a - b - c$ both lines are weak. Other examples of strong/weak arrows/lines are given in Figures 3, 4 and 7. (Also see [5, 6].)

The following fact will be used repeatedly. Consider the four statements:

1. $a \to b \in \mathcal{G}^*$;
2. $a \Rightarrow b \in G$ for all $G \in \mathcal{G}$;
3. $a \to b \in G$ for some $G \in \mathcal{G}$;
4. $a \Rightarrow b \in \mathcal{G}^*$.

Then 1 implies 2 implies 4 and, equivalently, 1 implies 3 implies 4.

This section will culminate with Theorem 3.2, which establishes that $\mathcal{G}^*$ is in fact adicyclic so is itself a chain graph. Thus, $\mathcal{G}^* \in \mathcal{G}$ by the following Lemma 3.1 and, therefore, $\mathcal{G}^*$ uniquely represents $\mathcal{G}$: $\mathcal{G}_1 = \mathcal{G}_2$ iff $\mathcal{G}_1^* = \mathcal{G}_2^*$.

LEMMA 3.1. *$\mathcal{G}^*$ has the same skeleton and triplexes as each $G \in \mathcal{G}$.*

PROOF. Clearly, $\mathcal{G}^*$ has the same skeleton as each $G \in \mathcal{G}$. Suppose that the triplex $(\{a, c\}, b)$ occurs in each $G \in \mathcal{G}$, so $a \Rightarrow b \Leftarrow c$ occurs as an induced subgraph in each $G \in \mathcal{G}$, hence, also $a \Rightarrow b \Leftarrow c$ occurs as an induced subgraph in $\mathcal{G}^*$. Furthermore, there exists $G \in \mathcal{G}$ such that either $a \to b \in G$, in which case $a \to b \in \mathcal{G}^*$, or $b \leftarrow c \in G$, in which case $b \leftarrow c \in \mathcal{G}^*$. In both cases, the triplex $(\{a,c\},b)$ occurs in $\mathcal{G}^*$, as required.

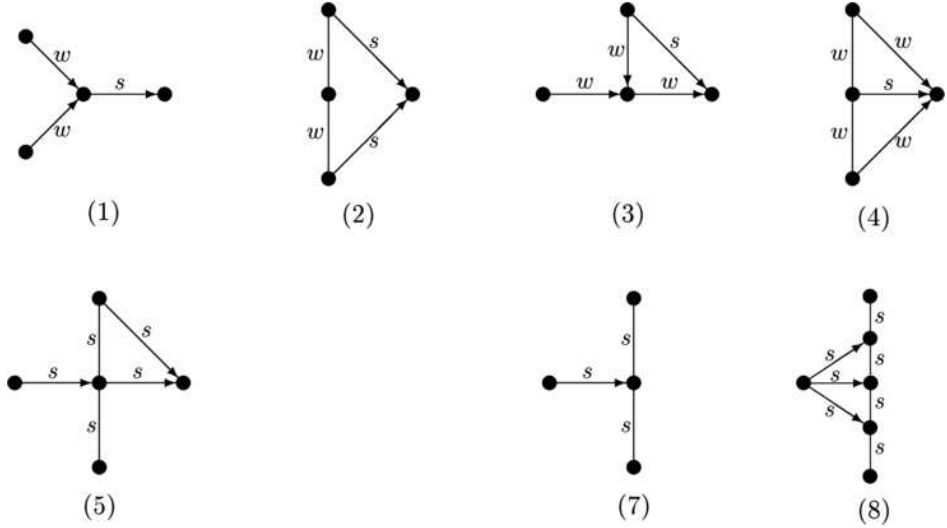

Fig. 4. *Strong/weak arrows/lines in AMP essential graphs (cf. Figure 6).*

Conversely, suppose that the triplex $(\{a,c\},b)$ occurs in $\mathcal{G}^*$. If this triplex occurs as the immorality $a \to b \leftarrow c$ in $\mathcal{G}^*$, then $a \Rightarrow b \Leftarrow c$ occurs as an induced subgraph in each $G \in \mathcal{G}$ while $a \to b \in G'$ for some $G' \in \mathcal{G}$, hence, this triplex occurs in $G'$ and thus in all $G \in \mathcal{G}$.

If this triplex occurs as the flag $a \to b - c$ in $\mathcal{G}^*$, then $a \Rightarrow b \in G$ for every $G \in \mathcal{G}$. If $b - c$ were a weak line in $\mathcal{G}^*$, then there would exist $G, G' \in \mathcal{G}$ such that $a \Rightarrow b \to c$ occurs in $G$ and $a \Rightarrow b \leftarrow c$ occurs in $G'$, both as induced subgraphs, contradicting the fact that $G$ and $G'$ have the same triplexes. Therefore, $b - c$ must be a strong line in $\mathcal{G}^*$. By a similar argument, $a \to b$ must be a strong arrow in $\mathcal{G}^*$ and, therefore, $a \to b - c$ occurs as a flag in each $G \in \mathcal{G}$. Thus, again the triplex $(\{a,c\},b)$ occurs in each $G \in \mathcal{G}$. □

LEMMA 3.2. (a) *A flag $a \to b - c$ occurs in $\mathcal{G}^*$ iff it occurs in each $G \in \mathcal{G}$, in which case $a \to b$ is a strong arrow and $b - c$ is a strong line. Thus, neither $a \to b \overset{\text{w}}{-} c$ nor $a \overset{\text{w}}{\to} b \overset{\text{s}}{-} c$ can occur as an induced subgraph of $\mathcal{G}^*$.*

(b)
- *If $a \to b \overset{\text{w}}{-} c$ occurs as a subgraph in $\mathcal{G}^*$, then $a \Rightarrow c \in \mathcal{G}^*$.*
- *If $a \overset{\text{w}}{\to} b \overset{\text{s}}{-} c$ occurs as a subgraph in $\mathcal{G}^*$, then $a \overset{\text{w}}{\to} c \in \mathcal{G}^*$.*
- *If $a \overset{\text{s}}{\to} b \overset{\text{w}}{\to} c$ occurs as a subgraph in $\mathcal{G}^*$, then $a \overset{\text{s}}{\to} c \in \mathcal{G}^*$.*
- *If $a \overset{\text{s}}{-} b \overset{\text{w}}{-} c$ occurs as a subgraph in $\mathcal{G}^*$, then $a \overset{\text{w}}{-} c \in \mathcal{G}^*$.*

PROOF. (a) "Only if" was established in the last paragraph of the proof of Lemma 3.1, while "if" is immediate.

(b) If $a \to b \overset{\text{w}}{-} c$ occurs in $\mathcal{G}^*$, then by (a) $a \cdots c \in \mathcal{G}^*$ and, therefore, $a \cdots c$ occurs in all $G \in \mathcal{G}$. Choose $G \in \mathcal{G}$ such that $b \to c \in G$. Since necessarily $a \Rightarrow b \in G$, therefore, $a \to c \in G$ by adicyclicity, hence $a \Rightarrow c \in \mathcal{G}^*$.

If $a \overset{\text{w}}{\to} b \overset{\text{s}}{-} c$ occurs in $\mathcal{G}^*$, then by (a) the triangle $a \Rightarrow b - c \cdots a$ occurs in all $G \in \mathcal{G}$, hence $a \Rightarrow c \in G$ by adicyclicity. If we now choose $G, G' \in \mathcal{G}$ such that $a \to b \in G$ and $a - b \in G'$, then $a \to c \in G$ and $a - c \in G'$ by adicyclicity, hence $a \overset{\text{w}}{\to} c \in \mathcal{G}^*$.

If $a \overset{\text{s}}{\to} b \overset{\text{w}}{\to} c$ occurs in $\mathcal{G}^*$, then $a \to b \Rightarrow c$ occurs in every $G \in \mathcal{G}$ and $a \to b - c$ occurs in some $G' \in \mathcal{G}$. Since $\mathcal{G}^*$ and $G'$ have the same triplexes by Lemma 3.1, $a \cdots c \in G'$, so the triangle $a \to b \Rightarrow c \cdots a$ occurs in every $G \in \mathcal{G}$. Thus, $a \to c$ must occur in every $G \in \mathcal{G}$ by adicyclicity, hence $a \overset{\text{s}}{\to} c \in \mathcal{G}^*$.

If $a \overset{\text{s}}{-} b \overset{\text{w}}{-} c$ occurs in $\mathcal{G}^*$, then there exist $G, G' \in \mathcal{G}$ such that $a - b \to c$ occurs in $G$ and $a - b \leftarrow c$ occurs in $G'$. Since $G$ and $G'$ have the same triplexes, necessarily $a \cdots c \in G, G'$, hence $a \to c \in G$ and $a \leftarrow c \in G'$ by adicyclicity. Therefore, $a \overset{\text{w}}{-} c \in \mathcal{G}^*$. □

LEMMA 3.3. *A biflag or chordless undirected cycle occurs in $\mathcal{G}^*$ iff it occurs in at least one $G \in \mathcal{G}$, in which case it occurs in all $G \in \mathcal{G}$. Thus, if an arrow occurs in a biflag in $\mathcal{G}^*$ or some $G \in \mathcal{G}$, then it must be strong, while if a line occurs in some biflag or chordless undirected cycle in $\mathcal{G}^*$ or in some $G \in \mathcal{G}$, then it must be strong.*

PROOF. If a biflag $[a; c_1, \ldots, c_k]$ $(k \geq 3)$ or $[a, b; c_1, \ldots, c_k]$ $(k \geq 2)$ occurs in $\mathcal{G}^*$, then, because $\overset{\text{s}}{-} \overset{\text{w}}{-}$ cannot occur as an induced subgraph in $\mathcal{G}^*$ [Lemma 3.2(b)], either all undirected edges $c_i - c_{i+1}$ $(1 \leq i \leq k)$ are strong or all are weak. Since the flag $a \overset{\text{s}}{\to} c_{k-1} \overset{\text{s}}{-} c_k$ occurs in $\mathcal{G}^*$, therefore all edges $c_i - c_{i+1}$ are strong. Furthermore, if any arrow $a \to c_i$ or $b \to c_i$ is weak, then all must be weak [apply Lemma 3.2(b)], a contradiction. Thus, all arrows and lines occurring in the biflag are strong, hence the biflag occurs in every $G \in \mathcal{G}$.

Conversely, if a biflag $[a; c_1, \ldots, c_k]$ or $[a, b; c_1, \ldots, c_k]$ occurs in some $G \in \mathcal{G}$, then the triplex $(\{a, c_k\}, c_{k-1})$ occurs in $\mathcal{G}^*$. If this triplex occurs as the flag $a \overset{\text{s}}{\to} c_{k-1} \overset{\text{s}}{-} c_k$, then the above argument applies to show that the biflag occurs in every $G \in \mathcal{G}$ and in $\mathcal{G}^*$. But the triplex cannot occur in either of the other two possible configurations: If it did, then $c_{k-1} \leftarrow c_k \in \mathcal{G}^*$, which would require that $c_1 \leftarrow c_2 \in \mathcal{G}^*$. (If $c_1 \Rightarrow c_2 \in \mathcal{G}^*$, then at least one triplex would occur in the path $c_1 \Rightarrow c_2 \cdots c_{k-1} \leftarrow c_k$ in $\mathcal{G}^*$, whereas no such triplex occurs in $G$, contradicting Lemma 3.1.) But $c_1 \leftarrow c_2 \notin \mathcal{G}^*$, since $a \to c_2 - c_1$ or $b \to c_2 - c_1$ occurs as a flag in $G$, so $(\{a, c_1\}, c_2) \in \mathcal{G}^*$ or $(\{b, c_1\}, c_2) \in \mathcal{G}^*$.

Next, suppose that a chordless undirected cycle occurs in $\mathcal{G}^*$ or in some $G \in \mathcal{G}$. Because no triplex occurs in this cycle, none can occur in the corresponding subgraph in any $G' \in \mathcal{G}$. Therefore, no arrow can occur in this subgraph in any $G' \in \mathcal{G}$; otherwise this subgraph would either include at least one triplex or else be a fully directed cycle, contradicting the adicyclicity of $G'$. Thus, the chordless undirected cycle occurs in every $G' \in \mathcal{G}^*$ and each of its lines is a strong line. $\square$

The converse to the second statement in Lemma 3.3 is not true: strong arrows not occurring in biflags appear in Figures 3(1) and 4(1)–(5). Strong lines not in biflags or chordless cycles may also occur; see [5, 6].

LEMMA 3.4. (a) *Any semi-directed cycle in $\mathcal{G}^*$ has at least one weak line.*
(b) *If $a \to b \Rightarrow d \Rightarrow a$ is a semi-directed 3-cycle in $\mathcal{G}^*$, then $\exists G, G' \in \mathcal{G}$ such that $a \Rightarrow b \leftarrow d \to a$ occurs in $G$ and $a \Rightarrow b \to d \leftarrow a$ occurs in $G'$. Thus, the semi-directed 3-cycle must have the form $a \to b \stackrel{\mathrm{w}}{-} d \stackrel{\mathrm{w}}{-} a$ in $\mathcal{G}^*$.*
(c) *In this case, $a \to b$ cannot occur in $\mathcal{G}^*$ in an immorality $a \to b \leftarrow c$, in a flag $a \to b - c$, in a chordless 2-dipath $c \to a \to b$ or in an antiflag of the form $c \stackrel{\mathrm{s}}{-} a \to b$.*

PROOF. (a) Let $d_0 \to d_1 \Rightarrow \cdots \Rightarrow (d_k \equiv d_0)$ ($k \geq 3$) be a semi-directed $k$-cycle in $\mathcal{G}^*$. Then $\exists G \in \mathcal{G}$ such that $d_0 \to d_1 \in G$. Since $G$ is adicyclic, $d_i \leftarrow d_{i+1} \in G$ for some $1 \leq i \leq k-1$, hence $d_i \Leftarrow d_{i+1} \in \mathcal{G}^*$. But $d_i \Rightarrow d_{i+1} \in \mathcal{G}^*$, so $d_i \stackrel{\mathrm{w}}{-} d_{i+1} \in \mathcal{G}^*$ as required.

(b) By (a), either $b \stackrel{\mathrm{w}}{-} d \in \mathcal{G}^*$ or $d \stackrel{\mathrm{w}}{-} a \in \mathcal{G}^*$ (or both). If the former, then $\exists G' \in \mathcal{G}$ such that $b \to d \in G'$. Since $a \Rightarrow b \in G'$, $d \leftarrow a \in G'$ by adicyclicity. Therefore, $d \Leftarrow a \in \mathcal{G}^*$, hence $d \stackrel{\mathrm{w}}{-} a \in \mathcal{G}^*$ (since $d \Rightarrow a \in \mathcal{G}^*$ and $d \leftarrow a \in G'$), so $\exists G \in \mathcal{G}$ such that $d \to a \in G$. Since $a \Rightarrow b \in G$, $b \leftarrow d \in G$ by adicyclicity. If the latter, then $\exists G \in \mathcal{G}$ such that $d \to a \in G$. Since $a \Rightarrow b \in G$, necessarily $b \leftarrow d \in G$ by adicyclicity. Thus, $b \stackrel{\mathrm{w}}{-} d \in \mathcal{G}^*$ (since $b \Rightarrow d \in \mathcal{G}^*$), so $\exists G' \in \mathcal{G}$ such that $b \to d \in G'$. Since $a \Rightarrow b \in G'$, $d \leftarrow a \in G'$ by adicyclicity.

(c) Assume that $a \to b \Leftarrow c$ occurs as an immorality or flag in $\mathcal{G}^*$. Let $G, G'$ be as specified in (b). Since $a \not\cdot c$, necessarily $c \neq d$. Because $\mathcal{G}^*, G, G'$ have the same triplexes (Lemma 3.1), $b \Leftarrow c \in G, G'$. Necessarily $d \cdots c$, for otherwise the triplex $(\{d,c\},b)$ would occur as the induced subgraph $d \to b \leftarrow c$ in $G$ but this triplex could not occur in $G'$. Thus, by adicyclicity, $d \leftarrow c \in G'$, so the triplex $(\{a,c\},d)$ occurs as the immorality $a \to d \leftarrow c$ in $G'$, but this triplex cannot occur in $G$, a contradiction.

Next, assume that either (i) the chordless 2-dipath $c \to a \to b$ or (ii) the antiflag $c \stackrel{\mathrm{s}}{-} a \to b$ occurs in $\mathcal{G}^*$. Let $G, G'$ be as specified in (b); again, necessarily $c \neq d$. For both (i) and (ii), necessarily $c \Rightarrow a \in G, G'$. Therefore, $c \cdots d$, for otherwise the triplex $(\{c,d\}, a)$ would occur as the induced subgraph $c \Rightarrow a \Leftarrow d$ in $G$, but this triplex could not occur in $G'$. Thus, by adicyclicity, $c \to d \in G'$, so the triplex $(\{c,b\}, d)$ occurs as the immorality $c \to d \leftarrow b$ in $G'$, but this triplex cannot occur in $G$, a contradiction. □

For any graph $H \equiv (W, F)$, define $H^\circ \equiv (W, F^\circ)$ to be the smallest chain graph larger than $H$, that is, $H^\circ$ is obtained from $H$ by converting any arrow that occurs in a semi-directed cycle in $H$ into a line. Note that this can be done in a single step: if, after converting an arrow that occurs in a semi-directed cycle in $H$ into a line, a second arrow now becomes part of a semi-directed cycle, then this second arrow already must have occurred in a semi-directed cycle in $H$. Clearly, $H^\circ$ has the same skeleton as $H$, $H \subseteq H^\circ$, and $H$ is adicyclic iff $H^\circ = H$.

LEMMA 3.5. (a) $\mathcal{G}^*$ and $(\mathcal{G}^*)^\circ$ have the same immoralities.
(b) $\mathcal{G}^*$ and $(\mathcal{G}^*)^\circ$ have the same flags.
(c) If an antiflag of the form $a \overset{\text{s}}{-} b \to c$ occurs $\mathcal{G}^*$, it also occurs in $(\mathcal{G}^*)^\circ$.

PROOF. (a) Since $\mathcal{G}^* \subseteq (\mathcal{G}^*)^\circ$ with the same skeletons, if $a \to b \leftarrow c$ occurs as an immorality in $(\mathcal{G}^*)^\circ$, then it must occur in $\mathcal{G}^*$. Conversely, suppose that $a \to b \leftarrow c$ occurs as an immorality in $\mathcal{G}^*$ but not in $(\mathcal{G}^*)^\circ$. This can happen only if at least one of the two arrows is converted to a line in $(\mathcal{G}^*)^\circ$, hence occurs in some semi-directed cycle in $\mathcal{G}^*$. *Choose the immorality $a \to b \leftarrow c$ that is associated with a semi-directed cycle of minimum length $k$* [i.e., minimal with respect to *all* immoralities $a' \to b' \leftarrow c'$ in $\mathcal{G}^*$ that do not occur in $(\mathcal{G}^*)^\circ$] and assume without loss of generality that this minimum-length semi-directed cycle contains $a \to b$. This cycle thus has the form $a \to (b \equiv d_1) \Rightarrow d_2 \Rightarrow \cdots \Rightarrow (d_k \equiv a)$ in $\mathcal{G}^*$. By Lemma 3.4(c), $k \geq 4$.

It is conceivable that $d_i = c$ for some (at most one) $i = 3, \ldots, k-2$. In that case, however, $(d_i \equiv c) \to b$ occurs in a shorter semi-directed cycle $c \to (b \equiv d_1) \Rightarrow \cdots \Rightarrow (d_i \equiv c)$ in $\mathcal{G}^*$, contradicting the minimality of $k$. Therefore, $d_i \neq c$ for each $i$.

By Lemma 3.4(a), the minimum-length semi-directed cycle has at least one weak line $d_j \overset{\text{w}}{-} d_{j+1} \in \mathcal{G}^*$ ($1 \leq j \leq k-1$). Consider the least such $j$.

Suppose first that $j \geq 2$. In this case the minimality of $j$ implies that either $d_{j-1} \to d_j \overset{\text{w}}{-} d_{j+1}$ or $d_{j-1} \overset{\text{s}}{-} d_j \overset{\text{w}}{-} d_{j+1}$ occurs in $\mathcal{G}^*$, hence $d_{j-1} \Rightarrow d_{j+1} \in \mathcal{G}^*$ by Lemma 3.2(b). Thus, $a \overset{\text{s}}{\to} (b \equiv d_1) \Rightarrow \cdots \Rightarrow d_{j-1} \Rightarrow d_{j+1} \Rightarrow \cdots \Rightarrow (d_k \equiv a)$ is a shorter semi-directed cycle in $\mathcal{G}^*$, contradicting the minimality of $k$.

Suppose now that $j = 1$, so $(b \equiv d_1) \overset{\text{w}}{-} d_2 \in \mathcal{G}^*$. Therefore, $a \to b \overset{\text{w}}{-} d_2$ and $d_2 \overset{\text{w}}{-} b \leftarrow c$ occur as subgraphs in $\mathcal{G}^*$, so by Lemma 3.2(b), $a \Rightarrow d_2 \Leftarrow c$ occurs as a subgraph of $\mathcal{G}^*$. If $a - d_2 \in \mathcal{G}^*$, then $a \to b$ occurs in the semi-directed 3-cycle $a \to b - d_2 - a$ in $\mathcal{G}^*$, contradicting Lemma 3.4(c), hence $a \to d_2 \in \mathcal{G}^*$.

Similarly, $d_2 \leftarrow c \in \mathcal{G}^*$, hence $a \to d_2 \leftarrow c$ occurs as an immorality in $\mathcal{G}^*$. But then $a \to d_2$ occurs in the shorter semi-directed cycle $a \to d_2 \Rightarrow \cdots \Rightarrow (d_k \equiv a)$ in $\mathcal{G}^*$, contradicting the minimality of $k$. We conclude that every immorality in $\mathcal{G}^*$ also occurs in $(\mathcal{G}^*)^\circ$.

(b) If a flag occurs in $\mathcal{G}^*$, it has the form $a \xrightarrow{\text{s}} b \xrightarrow{\text{s}} c$ by Lemma 3.2(a). Because $\mathcal{G}^* \subseteq (\mathcal{G}^*)^\circ$ and they have the same skeletons, if $a \to b - c$ does not occur as a flag in $(\mathcal{G}^*)^\circ$, then $a - b - c$ must occur in $(\mathcal{G}^*)^\circ$ with $a \cdot\!/\!\cdot c$ in both graphs. This implies that $a \xrightarrow{\text{s}} b$ occurs in a semi-directed cycle $a \xrightarrow{\text{s}} (b \equiv d_1) \Rightarrow d_2 \Rightarrow \cdots \Rightarrow (d_k \equiv a)$ in $\mathcal{G}^*$. Choose this cycle to minimize $k$; by Lemma 3.4(c), $k \geq 4$.

It is conceivable that $d_i = c$ for some (at most one) $i = 2, \ldots, k - 2$. If such $i$ exists and $i \geq 3$, then $a \xrightarrow{\text{s}} (b \equiv d_1) \xrightarrow{\text{s}} d_i \Rightarrow \cdots \Rightarrow (d_k \equiv a)$ is a shorter semi-directed cycle in $\mathcal{G}^*$, contradicting the minimality of $k$. Thus, either $i = 2$ or no such $i$ exists.

By Lemma 3.4(a), the minimum-length semi-directed cycle must have at least one weak line $d_j \xrightarrow{\text{w}} d_{j+1} \in \mathcal{G}^*$ ($1 \leq j \leq k - 1$). Consider the least such $j$. If $j \geq 2$, then the minimality of $k$ is contradicted exactly as in (a). Suppose, therefore, that $j = 1$, so $(b \equiv d_1) \xrightarrow{\text{w}} d_2 \in \mathcal{G}^*$. (Also, $i$ does not exist and $d_2 \neq c$.) Since now $a \xrightarrow{\text{s}} b \xrightarrow{\text{w}} d_2$ and $d_2 \xrightarrow{\text{w}} b \xrightarrow{\text{s}} c$ occur in $\mathcal{G}^*$, $a \Rightarrow d_2 \in \mathcal{G}^*$ and $d_2 \xrightarrow{\text{w}} c \in \mathcal{G}^*$ by Lemma 3.2(b). However, $a \to d_2$ cannot occur in $\mathcal{G}^*$ [otherwise $a \to d_2 \xrightarrow{\text{w}} c$ would be a flag in $\mathcal{G}^*$, contradicting Lemma 3.2(a)], so $a - d_2 \in \mathcal{G}^*$. Thus, $a \to b - d_2 - a$ is a semi-directed 3-cycle in $\mathcal{G}^*$, contrary to Lemma 3.4.

Conversely, assume that $a \to b - c$ occurs as a flag in $(\mathcal{G}^*)^\circ$, so $a \to b \cdots c$ occurs in $\mathcal{G}^*$. If $b - c \in \mathcal{G}^*$, then the flag also occurs in $\mathcal{G}^*$, as asserted. If $b \leftarrow c \in \mathcal{G}^*$, then $a \to b \leftarrow c$ occurs as an immorality in $\mathcal{G}^*$, hence by part (a) also occurs as an immorality in $(\mathcal{G}^*)^\circ$, contrary to assumption.

Thus, assume that $a \to b \to c$ occurs as a chordless 2-dipath in $\mathcal{G}^*$. Then $b \to c$ must occur in a semi-directed cycle $b \to (c \equiv d_1) \Rightarrow d_2 \Rightarrow \cdots \Rightarrow (d_k \equiv b)$ in $\mathcal{G}^*$. [Note that $d_i \neq a$ for each $i$; otherwise $a \to b$ would be included in a semi-directed cycle and would thus be converted to $a - b$ in $(\mathcal{G}^*)^\circ$, contradicting the original assumption.] *Choose the chordless 2-dipath $a \to b \to c$ in $\mathcal{G}^*$ that is associated with a semi-directed cycle having minimal length $k$ with respect to all chordless 2-dipaths $a' \to b' \to c'$ in $\mathcal{G}^*$ that occur as $a' \to b' - c$ in $(\mathcal{G}^*)^\circ$.* By Lemma 3.4(c), $k \geq 4$.

By Lemma 3.4(a), the minimum-length semi-directed cycle must have at least one weak line $d_j \xrightarrow{\text{w}} d_{j+1} \in \mathcal{G}^*$ ($1 \leq j \leq k - 1$); consider the least such $j$. If $j \geq 2$, then the minimality of $k$ is again contradicted exactly as in (a). Thus, assume that $j = 1$, so $b \to (c \equiv d_1) \xrightarrow{\text{w}} d_2$ occurs in $\mathcal{G}^*$. By Lemma 3.2(b), $b \Rightarrow d_2 \in \mathcal{G}^*$. If $b - d_2 \in \mathcal{G}^*$, then $b \to c$ occurs in a semi-directed 3-cycle $b \to c \xrightarrow{\text{w}} d_2 - b$ in $\mathcal{G}^*$, contradicting Lemma 3.4(c), so $b \to d_2 \in \mathcal{G}^*$.

If $a \cdot\!/\!\cdot d_2$ in $\mathcal{G}^*$, then the chordless 2-dipath $a \to b \to d_2$ is associated with a shorter semi-directed cycle $b \to d_2 \Rightarrow \cdots \Rightarrow (d_k \equiv b)$ in $\mathcal{G}^*$, contradicting the minimality of $k$, hence $a \cdots d_2 \in \mathcal{G}^*$. If $a \Leftarrow d_2 \in \mathcal{G}^*$, then $a \to b \to c \stackrel{\text{w}}{-} d_2 \Rightarrow a$ is a semi-directed cycle in $\mathcal{G}^*$, so $\alpha \to \beta$ would be converted to a line in $(\mathcal{G}^*)^\circ$, contrary to assumption; thus, $a \to d_2 \in \mathcal{G}^*$. But then $a \to d_2 \stackrel{\text{w}}{-} c$ is a flag in $\mathcal{G}^*$, contradicting Lemma 3.2(a). Thus, $a \to b - c$ also occurs as a flag in $\mathcal{G}^*$.

(c) Suppose that $a \stackrel{\text{s}}{-} b \to c$ occurs as an antiflag in $\mathcal{G}^*$ but not in $(\mathcal{G}^*)^\circ$. Then $b \to c$ must occur in a semi-directed cycle $b \to (c \equiv d_1) \Rightarrow d_2 \Rightarrow \cdots \Rightarrow (d_k \equiv b)$ in $\mathcal{G}^*$. Choose the antiflag $a \stackrel{\text{s}}{-} b \to c$ in $\mathcal{G}^*$ that is associated with a semi-directed cycle having minimum length $k$ with respect to all antiflags $a' \stackrel{\text{s}}{-} b' \to c'$ in $\mathcal{G}^*$ that occur as $a' - b' - c$ in $(\mathcal{G}^*)^\circ$. By Lemma 3.4(c), $k \geq 4$.

It is conceivable that $d_i = a$ for some (at most one) $i = 3, \ldots, k-1$. If such $i$ exists and $i \leq k-2$, then $b \to (c \equiv d_1) \Rightarrow \cdots \Rightarrow d_i \stackrel{\text{s}}{-} (d_k \equiv b)$ is a shorter semi-directed cycle in $\mathcal{G}^*$, contradicting the minimality of $k$. Thus, either $i = k - 1$ or no such $i$ exists.

By Lemma 3.4(a), the minimum-length semi-directed cycle must have at least one weak line $d_j \stackrel{\text{w}}{-} d_{j+1} \in \mathcal{G}^*$ ($1 \leq j \leq k-1$); consider the least such $j$. If $j \geq 2$, then the minimality of $k$ is again contradicted exactly as in (a). Thus, assume that $j = 1$, so $b \to (c \equiv d_1) \stackrel{\text{w}}{-} d_2$ occurs in $\mathcal{G}^*$. By Lemma 3.2(b), $b \Rightarrow d_2 \in \mathcal{G}^*$. If $b - d_2 \in \mathcal{G}^*$, then $b \to c$ occurs in a semi-directed 3-cycle $b \to c \stackrel{\text{w}}{-} d_2 - b$ in $\mathcal{G}^*$, contradicting Lemma 3.4(c), so $b \to d_2 \in \mathcal{G}^*$.

If $a \cdot / \cdot d_2$ in $\mathcal{G}^*$, then the antiflag $a \stackrel{\text{s}}{\to} b \to d_2$ is associated with a shorter semi-directed cycle $b \to d_2 \Rightarrow \cdots \Rightarrow (d_k \equiv b)$ in $\mathcal{G}^*$, contradicting the minimality of $k$; hence $a \cdots d_2 \in \mathcal{G}^*$. If $a \Leftarrow d_2 \in \mathcal{G}^*$, then $b \to d_2 \Rightarrow a \stackrel{\text{s}}{-} b$ is a semi-directed 3-cycle in $\mathcal{G}^*$, contradicting Lemma 3.4(b), so $a \to d_2 \in \mathcal{G}^*$. But then $a \to d_2 \stackrel{\text{w}}{-} c$ occurs as a flag in $\mathcal{G}^*$, contradicting Lemma 3.2(a). □

Two vertices $a, a' \in V$ are *strongly equivalent* (with respect to $\mathcal{G}$) if $a = a'$ or there is a path between them in $\mathcal{G}^*$ consisting solely of strong lines. Let $\Sigma \equiv \Sigma(\mathcal{G})$ denote the set of strong equivalence classes in $\mathcal{G}^*$, providing the decomposition $V = \dot{\bigcup}(\sigma | \sigma \in \Sigma)$. If $a \cdots a' \in \mathcal{G}^*_\sigma$, then $a - a' \in G_\sigma \ \forall G \in \mathcal{G}$ by the adicyclicity of $G$, so $a \stackrel{\text{s}}{-} a' \in \mathcal{G}^*_\sigma$. Therefore, $\mathcal{G}^*_\sigma = ((\mathcal{G}^*)^\circ)_\sigma = G_\sigma$ is a connected UG (possibly a singleton), each of whose lines is strong. For $a \in V$, the unique strong equivalence class containing $a$ is denoted by $\sigma(a)$.

LEMMA 3.6. *Suppose that $a \in \alpha$ and $b \in \beta$ for distinct $\alpha, \beta \in \Sigma \equiv \Sigma(\mathcal{G})$.*

(a) *If $a \to b \in G$ (resp., $a - b \in G$) for some $G \in \mathcal{G}$ and $a' \cdots b' \in G$ for a pair $a' \in \alpha$, $b' \in \beta$, then $a' \to b' \in G$ (resp., $a' - b' \in G$).*

(b) *If $a \stackrel{\text{w}}{\to} b \in \mathcal{G}^*$ (resp., $a \stackrel{\text{s}}{\to} b \in \mathcal{G}^*$) and $a' \cdots b' \in \mathcal{G}^*$ for a pair $a' \in \alpha$, $b' \in \beta$, then $a' \stackrel{\text{w}}{\to} b' \in \mathcal{G}^*$ (resp., $a' \stackrel{\text{s}}{\to} b' \in \mathcal{G}^*$).*

(c) *If $a \stackrel{\text{w}}{\to} b \in \mathcal{G}^*$, then $a \stackrel{\text{w}}{\to} b' \in \mathcal{G}^*$ for every $b' \in \beta$. Furthermore, $\mathcal{G}^*_{\alpha \cap \mathrm{wp}_{\mathcal{G}^*}(\beta)}$ is complete.*

(d) If $a \stackrel{\text{w}}{-} b \in \mathcal{G}^*$, then $a' \stackrel{\text{w}}{-} b' \in \mathcal{G}^*$ for every pair $a' \in \alpha$, $b' \in \beta$. Furthermore, $\mathcal{G}^*_\alpha$ and $\mathcal{G}^*_\beta$ are complete.

PROOF. (a) is immediate by the adicyclicity of $G$ and the connectedness of the subgraphs $G_\alpha$ and $G_\beta$.

(b) If $a \stackrel{\text{w}}{\to} b \in \mathcal{G}$, then $a \to b \in G_1$ and $a - b \in G_2$ for some $G_1, G_2 \in \mathcal{G}$, while $a \Rightarrow b \in G_3$ for all other $G_3 \in \mathcal{G}$. Therefore, by (a), $a' \to b' \in G_1$, $a' - b' \in G_2$, and $a' \Rightarrow b' \in G_3$, hence $a' \stackrel{\text{w}}{\to} b' \in \mathcal{G}^*$. Similarly, if $a \stackrel{\text{s}}{\to} b \in \mathcal{G}^*$, then $a' \stackrel{\text{s}}{\to} b' \in \mathcal{G}^*$.

(c) By Lemma 3.2(b), $a \stackrel{\text{w}}{\to} b'' \in \mathcal{G}^*$ for all strong neighbors $b''$ of $b$ in $\mathcal{G}^*$, so by the connectedness of $\mathcal{G}^*_\beta$, $a \stackrel{\text{w}}{\to} b' \in \mathcal{G}^*$ for every $b' \in \beta$. Next, suppose that $a', a'' \in \alpha \cap \text{wp}_{\mathcal{G}^*}(\beta)$. It follows from the preceding that $a' \stackrel{\text{w}}{\to} b \in \mathcal{G}^*$ and $a'' \stackrel{\text{w}}{\to} b \in \mathcal{G}^*$. By (a), therefore, $\exists G_1, G_2 \in \mathcal{G}$ such that $a' \to b \in G_1$ and $a'' \to b \in G_1$, but $a' - b \in G_2$ and $a'' - b \in G_2$. Thus, $a' \cdots a'' \in \mathcal{G}^*$, for otherwise the triplex $(\{a', a''\}, b)$ would occur as the immorality $a' \to b \leftarrow a''$ in $G_1$, but this triplex would not occur in $G_2$.

(d) By Lemma 3.2(b), $a \stackrel{\text{w}}{-} b'' \in \mathcal{G}^*$ for all strong neighbors $b''$ of $b$ in $\mathcal{G}^*$, so by the connectedness of $G_\beta$, $a \stackrel{\text{w}}{-} b' \in \mathcal{G}^*$ for every $b' \in \beta$. Similarly, $a' \stackrel{\text{w}}{-} b' \in \mathcal{G}^*$ for all $a' \in \alpha$. Next, for any $a', a'' \in \mathcal{G}^*_\alpha$, the preceding shows that $a' \stackrel{\text{w}}{-} b \in \mathcal{G}^*$ and $a'' \stackrel{\text{w}}{-} b \in \mathcal{G}^*$. By (a), $\exists G_1, G_2 \in \mathcal{G}$ such that $a' \to b \in G_1$ and $a'' \to b \in G_1$, but $a' \leftarrow b \in G_2$ and $a'' \leftarrow b \in G_2$. Thus, $a' \cdots a'' \in \mathcal{G}^*$, for otherwise the triplex $(\{a', a''\}, b)$ would occur as the immorality $a' \to b \leftarrow a''$ in $G_1$, but this triplex would not occur in $G_2$, hence $\mathcal{G}^*_\alpha$ is complete. Similarly, $\mathcal{G}^*_\beta$ is complete. $\square$

DEFINITION 3.2. For each $G \equiv (V, E) \in \mathcal{G}$, define the *reduced graph* $\overline{G} \equiv (\Sigma \equiv \Sigma(\mathcal{G}), \overline{E})$ as follows: $\alpha \to \beta \in \overline{G}$ (resp., $\alpha - \beta \in \overline{G}$) iff $a \to b \in G$ (resp., $a - b \in G$) for at least one pair $a \in \alpha$, $b \in \beta$. The reduced graph $\overline{\mathcal{G}^*} \equiv (\Sigma, \overline{E^*})$ is defined similarly.

By Lemma 3.6(a), $\overline{G}$ is well defined. Since $G$ and $\overline{G}$ are adicyclic, $\overline{G^\circ} = \overline{G} = (\overline{G})^\circ$. Clearly, $\overline{\mathcal{G}^*}$ has the same skeleton as $\overline{G}$ for each $G \in \mathcal{G}$, and $\alpha \cdots \beta \in \overline{\mathcal{G}^*}$ iff $a \cdots b \in \mathcal{G}^*$ for some pair $a \in \alpha$, $b \in \beta$. The arrow $\alpha \to \beta \in \overline{\mathcal{G}^*}$ is called *weak* (resp., *strong*) if $a \stackrel{\text{w}}{\to} b \in \mathcal{G}^*$ (resp., $a \stackrel{\text{s}}{\to} b \in \mathcal{G}^*$). All lines $\alpha - \beta \in \overline{\mathcal{G}^*}$ are called *weak* because each line $a - b \in \mathcal{G}^*$ ($a \in \alpha$, $b \in \beta$) must be weak.

The following fact will be used repeatedly. Consider the following four statements:

1. $\alpha \to \beta \in \overline{\mathcal{G}^*}$;
2. $\alpha \Rightarrow \beta \in \overline{G}$ for all $G \in \mathcal{G}$;
3. $\alpha \to \beta \in \overline{G}$ for some $G \in \mathcal{G}$;

4. $\alpha \Rightarrow \beta \in \overline{\mathcal{G}^*}$.

Then $1 \Longleftrightarrow 2 \Longrightarrow 3 \Longleftrightarrow 4$.

LEMMA 3.7. (a) *If $G \in \mathcal{G}$, then $\overline{G}$ has no chordless undirected cycles.*
(b) *If $\alpha \to \beta \leftarrow \gamma$ occurs as an immorality in $\overline{G}$ for some $G \in \mathcal{G}$, then the triplex $(\{\alpha, \gamma\}, \beta)$ occurs in $\overline{G'}$ for every $G' \in \mathcal{G}$.*

PROOF. (a) If $(\sigma_0, \sigma_1, \ldots, \sigma_k \equiv \sigma_0)$ comprises a chordless undirected $k$-cycle in $\overline{G}$, then, by the definition of $\overline{G}$ and the connectedness of each $\mathcal{G}_{\sigma_i}$, there exist $\{s_{ij}|1 \leq j \leq n_i\} \subseteq \sigma_i$, $i = 1, \ldots, k$, such that $(s_{11}, \ldots, s_{1n_1}, \ldots, s_{k1}, \ldots, s_{kn_k} \equiv s_{11})$ is an undirected $l$-cycle in $G$ ($l \geq k$). If we choose the $l$-cycle of this form that minimizes $l$, this cycle must be chordless. By Lemma 3.3, each line $s_{in_i} - s_{(i+1)1}$ must be strong in $\mathcal{G}^*$, but also must be weak in $\mathcal{G}^*$ since $s_{in_i} \in \sigma_i$ and $s_{(i+1)1} \in \sigma_{i+1}$, a contradiction.

(b) Because $\alpha \to \beta \leftarrow \gamma$ occurs as an immorality in $\overline{G}$, $\alpha \Rightarrow \beta \Leftarrow \gamma$ occurs in $\overline{\mathcal{G}^*}$ with $\alpha \cdot / \cdot \gamma$. If $\alpha \overset{w}{\to} \beta$, $\alpha \overset{w}{-} \beta$, $\beta \overset{w}{\leftarrow} \gamma$, or $\beta \overset{w}{-} \gamma$ occurs in $\overline{\mathcal{G}^*}$, then, by Lemma 3.6(a), (b), (c), (d), $\exists a \in \alpha$, $b \in \beta$, and $c \in \gamma$ such that $a \to b \leftarrow c$ occurs as an immorality in $G$. In this case $(\{a, c\}, b)$ must occur as a triplex in each $G' \in \mathcal{G}$, hence $(\{\alpha, \gamma\}, \beta)$ must occur as a triplex in $\overline{G'}$ for each $G' \in \mathcal{G}$.

In the remaining case, $\alpha \overset{s}{\to} \beta \overset{s}{\leftarrow} \gamma$ occurs in $\mathcal{G}^*$. By the connectedness of $\mathcal{G}^*_\beta$, $\exists a \in \alpha$, $b_0, \ldots, b_n \in \beta$ ($n \geq 0$), and $c \in \gamma$ such that $a \overset{s}{\to} b_0 \overset{s}{-} \cdots \overset{s}{-} b_n \overset{s}{\leftarrow} c$ occurs as a subgraph of $\mathcal{G}^*$ and thus as a subgraph of each $G' \in \mathcal{G}$. Therefore, $\alpha \to \beta \leftarrow \gamma$ occurs as an immorality (a triplex) in $\overline{G'}$ for each $G' \in \mathcal{G}$. □

LEMMA 3.8. (a) *Let $\alpha, \beta \in \Sigma \equiv \Sigma(\mathcal{G})$ be distinct strong equivalence classes. An arrow $\alpha \to \beta \in \overline{\mathcal{G}^*}$ is contained in some semi-directed cycle in $\overline{\mathcal{G}^*}$ if and only if there exist $a \in \alpha$ and $b \in \beta$ such that $a \to b \in \mathcal{G}^*$ is contained in some semi-directed cycle in $\mathcal{G}^*$. In this case, every arrow $a' \to b' \in \mathcal{G}^*$ between $\alpha$ and $\beta$ is contained in some semi-directed cycle in $\mathcal{G}^*$.*
(b) *$\mathcal{G}^*$ is adicyclic iff $\overline{\mathcal{G}^*}$ is adicyclic. In general, $(\overline{\mathcal{G}^*})^\circ = \overline{(\mathcal{G}^*)^\circ}$.*

PROOF. (a) "Only if": This follows immediately from the fact that any two vertices in a strong equivalence class are connected in $\mathcal{G}^*$ via a path consisting of strong lines.

"If": Let $a \to b \Rightarrow d_1 \Rightarrow \cdots \Rightarrow (d_l \equiv a)$ be a semi-directed $l$-cycle in $\mathcal{G}^*$ ($l \geq 3$). Consider the largest $i = 1, \ldots, l-1$ such that $d_i \in \beta$. Since $d_i \Rightarrow d_{i+1} \in \mathcal{G}^*$, Lemma 3.6(b) implies that $\sigma(d_{i+1}) \neq \alpha$, hence $i \leq l - 2$, $\gamma := \sigma(d_{i+1}) \neq \alpha, \beta$, and $\alpha \to \beta \Rightarrow \gamma$ occurs as a subgraph in $\overline{\mathcal{G}^*}$. Next, consider the largest $j = i+1, \ldots, l-1$ such that $d_j \in \gamma$. Then $\delta := \sigma(d_{j+1}) \neq \beta, \gamma$ and, since $d_j \Rightarrow d_{j+1} \in \mathcal{G}^*$, $\alpha \to \beta \Rightarrow \gamma \Rightarrow \delta$ occurs in $\overline{\mathcal{G}^*}$. Either $\delta = \alpha$, producing a semi-directed 3-cycle of the desired form, or this process may be continued until a semi-directed $k$-cycle $\alpha \to \beta \Rightarrow \gamma \Rightarrow \delta \Rightarrow \cdots \Rightarrow \alpha$ $(k \geq 4)$ in $\overline{\mathcal{G}^*}$ is obtained. The final statement is immediate.

(b) The first part follows from the first statement in (a). By Lemma 3.5, $(\mathcal{G}^*)^\circ \in \mathcal{G}$, so $\overline{(\mathcal{G}^*)^\circ}$ is well defined; clearly, it has the same skeleton as $\overline{\mathcal{G}^*}$ and $(\overline{\mathcal{G}^*})^\circ$. Suppose first that $\alpha \to \beta \in (\overline{\mathcal{G}^*})^\circ$. Then $\alpha \to \beta \in \overline{\mathcal{G}^*}$ and $\alpha \to \beta$ cannot occur in a semi-directed cycle in $\overline{\mathcal{G}^*}$. Select $a \in \alpha$ and $b \in \beta$ such that $a \to b \in \mathcal{G}^*$. By (a)("if"), $a \to b$ cannot occur in a semi-directed cycle in $\mathcal{G}^*$, hence $a \to b \in (\mathcal{G}^*)^\circ$, so $\alpha \to \beta \in \overline{(\mathcal{G}^*)^\circ}$.

Suppose next that $\alpha - \beta \in (\overline{\mathcal{G}^*})^\circ$, so either (i) $\alpha - \beta \in \overline{\mathcal{G}^*}$, (ii) $\alpha \to \beta \in \overline{\mathcal{G}^*}$, or (iii) $\alpha \leftarrow \beta \in \overline{\mathcal{G}^*}$. In case (i), it follows directly that $\alpha - \beta \in \overline{(\mathcal{G}^*)^\circ}$. In case (ii), $\alpha \to \beta$ must occur in some semi-directed cycle in $\overline{\mathcal{G}^*}$, hence by (a)("only if") again $\alpha - \beta \in \overline{(\mathcal{G}^*)^\circ}$; case (iii) is similar to (ii). Thus, $(\overline{\mathcal{G}^*})^\circ = \overline{(\mathcal{G}^*)^\circ}$. □

LEMMA 3.9. (a) $\overline{\mathcal{G}^*}$ has no flags or chordless undirected cycles.
(b) If $\alpha \to \beta - \gamma$ occurs as a subgraph in $\overline{\mathcal{G}^*}$, then $\alpha \Rightarrow \gamma \in \overline{\mathcal{G}^*}$.
(c) Any triplex (necessarily an immorality) in $\overline{\mathcal{G}^*}$ is a triplex in $\overline{G}$ for each $G \in \mathcal{G}$.

PROOF. (a) If $\alpha \to \beta - \gamma$ is a flag in $\overline{\mathcal{G}^*}$, then $\beta - \gamma$ is a weak line and, by Lemma 3.6(d) and the definition of $\overline{\mathcal{G}^*}$, there exist $a \in \alpha$, $b \in \beta$, $c \in \gamma$ such that $a \to b - c$ is a flag in $\mathcal{G}^*$. Therefore, $b - c$ must be a strong line in $\mathcal{G}^*$ by Lemma 3.2(a), but $b - c$ is a weak line, thus, a contradiction.

Next, if $(\sigma_0, \sigma_1, \ldots, \sigma_k \equiv \gamma_0)$ is a chordless undirected $k$-cycle in $\overline{\mathcal{G}^*}$, then by Lemma 3.6(d) and the definition of $\overline{\mathcal{G}^*}$, there exist $s_i \in \sigma_i$, $i = 1, \ldots, k$, such that $(s_0, s_1, \ldots, s_k \equiv s_0)$ is a chordless undirected $k$-cycle in $\mathcal{G}^*$. Each line $s_{i-1} - s_i \in \mathcal{G}^*$ is weak, but, by Lemma 3.3, each line $s_{i-1} - s_i \in \mathcal{G}^*$ must be strong, a contradiction.

(b) By (a), $\alpha \cdots \gamma \in \overline{\mathcal{G}^*}$. Because $\beta - \gamma$ must be weak, we can choose $G \in \mathcal{G}$ such that $\beta \to \gamma \in \overline{G}$, hence, the triangle $\alpha \Rightarrow \beta \to \gamma \cdots \alpha$ occurs in $\overline{G}$. Since $\overline{G}$ is adicyclic, $\alpha \to \gamma \in \overline{G}$, hence $\alpha \Rightarrow \gamma \in \overline{\mathcal{G}^*}$.

(c) By (a), any triplex in $\overline{\mathcal{G}^*}$ is an immorality, say $\alpha \to \beta \leftarrow \gamma$. If either arrow is weak, then by Lemma 3.6(c) there exist $a \in \alpha$, $b \in \beta$, $c \in \gamma$ such that $a \to b \leftarrow c$ is an immorality in $\mathcal{G}^*$. Thus, by Lemma 3.1, $(\{a,c\}, b)$ is a triplex in $G$, so by the definition of $\overline{G}$, $(\{\alpha, \gamma\}, \beta)$ is a triplex in $\overline{G}$. If both arrows in $\alpha \to \beta \leftarrow \gamma$ are strong, then there exist $a \in \alpha$, $b, b' \in \beta$, $c \in \gamma$ such that $a \to b \in G$ and $b' \leftarrow c \in G$, hence $\alpha \to \beta \leftarrow \gamma$ is an immorality in $\overline{G}$. □

LEMMA 3.10. (a) Any semi-directed cycle in $\overline{\mathcal{G}^*}$ has at least one (weak) line.
(b) If $\alpha \to \beta \Rightarrow \delta \Rightarrow \alpha$ is a semi-directed 3-cycle in $\overline{\mathcal{G}^*}$, then there exist $G, G' \in \mathcal{G}$ such that $\alpha \Rightarrow \beta \leftarrow \delta \to \alpha$ occurs in $\overline{G}$ and $\alpha \Rightarrow \beta \to \delta \leftarrow \alpha$ occurs in $\overline{G'}$. Therefore, the semi-directed 3-cycle in $\overline{\mathcal{G}^*}$ must have the form $\alpha \to \beta \stackrel{w}{-} \delta \stackrel{w}{-} \alpha$.

(c) In this case, $\alpha \to \beta$ cannot occur in $\overline{\mathcal{G}^*}$ in a triplex (necessarily an immorality $\alpha \to \beta \leftarrow \gamma$) or in an induced subgraph of the form $\gamma \to \alpha \to \beta$.

(d) If $\sigma_0 \to \sigma_1 \Rightarrow \cdots \Rightarrow (\sigma_k \equiv \sigma_0)$ is a semi-directed $k$-cycle in $\overline{\mathcal{G}^*}$ ($k \geq 3$), then each $\sigma_i$ has at least one weak neighbor in $\mathcal{G}^*$. Thus, by Lemma 3.6(d), the induced subgraphs $\mathcal{G}^*_{\sigma_1}, \ldots, \mathcal{G}^*_{\sigma_k}$ are complete.

PROOF. (a) This follows directly from Lemma 3.4(a), using the fact that each $\mathcal{G}^*_\sigma$, $\sigma \in \Sigma$, is connected.

(b) The proof is the same as that of Lemma 3.4(b).

(c) Assume that the immorality $\alpha \to \beta \leftarrow \gamma$ occurs in $\overline{\mathcal{G}^*}$. Let $G, G'$ be as specified in (b); since $\alpha \cdot/\cdot \gamma$, necessarily $\gamma \neq \delta$. Because $\beta \leftarrow \gamma \in \overline{\mathcal{G}^*}$, necessarily $\beta \Leftarrow \gamma \in \overline{G}, \overline{G'}$. Since $\delta - \beta \leftarrow \gamma$ occurs as a subgraph of $\overline{\mathcal{G}^*}$, necessarily $\delta \Leftarrow \gamma \in \overline{\mathcal{G}^*}$ [Lemma 3.9(b)], hence also $\delta \cdots \gamma \in \overline{G}, \overline{G'}$. Thus, the triangle $\delta \leftarrow \beta \Leftarrow \gamma \cdots \delta$ occurs in $\overline{G'}$, hence $\delta \leftarrow \gamma \in \overline{G}$ by the adicyclicity of $\overline{G}$. Therefore, $\alpha \to \delta \leftarrow \gamma$ occurs as an immorality in $\overline{G'}$, so by Lemma 3.7(b), $(\{\alpha, \gamma\}, \delta)$ must occur as a triplex in $\overline{G}$. But this is impossible because $\alpha \leftarrow \delta \in \overline{G}$. The impossibility of the occurrence of $\gamma \to \alpha \to \beta$ as an induced subgraph of $\overline{\mathcal{G}^*}$ is proved by a similar argument.

(d) Use induction on $k$. By (b), the result is true when $k = 3$. Suppose it is true for all $k' < k$. By (a), at least one edge in the $k$-cycle is a weak line. Let $\sigma_j \overset{\text{w}}{-} \sigma_{j+1}$ be the first such edge ($1 \leq j \leq k-1$), so $\sigma_{j-1} \to \sigma_j \overset{\text{w}}{-} \sigma_{j+1}$ occurs in $\overline{\mathcal{G}^*}$. Thus, $\sigma_{j-1} \Rightarrow \sigma_{j+1} \in \overline{\mathcal{G}^*}$ by Lemma 3.9(b). If $j \geq 2$, then $\sigma_0 \to \sigma_1 \Rightarrow \cdots \Rightarrow \sigma_{j-1} \Rightarrow \sigma_{j+1} \Rightarrow \cdots \Rightarrow (\sigma_k \equiv \sigma_0)$ is a semi-directed $(k-1)$-cycle in $\overline{\mathcal{G}^*}$, so by the induction hypothesis, each of its vertices has at least one weak neighbor in $\overline{\mathcal{G}^*}$. Because the vertex $\sigma_j$ also has a weak neighbor ($\sigma_{j+1}$), the asserted result holds. If $j = 1$, then the cycle $\sigma_0 \Rightarrow \sigma_2 \Rightarrow \cdots \Rightarrow (\sigma_k \equiv \sigma_0)$ occurs in $\overline{\mathcal{G}^*}$. If $\sigma_0 \to \sigma_2 \in \overline{\mathcal{G}^*}$, then this cycle is semi-direct and the induction hypothesis again yields the asserted result. If $\sigma_0 \overset{\text{w}}{-} \sigma_2 \in \overline{\mathcal{G}^*}$, then either the cycle is completely undirected, in which case the asserted result obviously holds, or else it has at least one arrow, in which case the cycle is semi-directed and the induction hypothesis again applies. □

LEMMA 3.11. (a) *An immorality $\alpha \to \beta \leftarrow \gamma$ occurs in $\overline{\mathcal{G}^*}$ iff it occurs in $(\overline{\mathcal{G}^*})^\circ$.*

(b) *Neither $\overline{\mathcal{G}^*}$ nor $(\overline{\mathcal{G}^*})^\circ$ has any flags.*

(c) *$(\overline{\mathcal{G}^*})^\circ$ has no chordless undirected cycles.*

PROOF. (a) Note that $\overline{\mathcal{G}^*} \subseteq (\overline{\mathcal{G}^*})^\circ$ and they have the same skeletons. Thus, if the immorality $\alpha \to \beta \leftarrow \gamma$ occurs in $(\overline{\mathcal{G}^*})^\circ$, it must occur in $\overline{\mathcal{G}^*}$. Conversely, suppose that $\alpha \to \beta \leftarrow \gamma$ occurs in $\overline{\mathcal{G}^*}$ but not in $(\overline{\mathcal{G}^*})^\circ$. Then at least one of these two arrows must occur in a semi-directed cycle in $\overline{\mathcal{G}^*}$ and so is converted to a line in $(\overline{\mathcal{G}^*})^\circ$; also $\alpha \cdot/\cdot \gamma$ in both graphs. Choose the immorality $\alpha \to \beta \leftarrow \gamma$ associated with a semi-directed cycle of *minimum length* $k$ [minimum with respect to *all* immoralities $\alpha' \to \beta' \leftarrow \gamma'$ in $\overline{\mathcal{G}^*}$ that do not occur in $(\overline{\mathcal{G}^*})^\circ$] and assume without loss of generality that this cycle contains $\alpha \to \beta$. This cycle thus has the form $\alpha \to (\beta \equiv \delta_1) \Rightarrow \delta_2 \Rightarrow \cdots \Rightarrow (\delta_k \equiv \alpha)$ in $\overline{\mathcal{G}^*}$. By Lemma 3.10(c), $k \geq 4$.

It may be that $\delta_i = \gamma$ for some (at most one) $i = 3, \ldots, k-2$. In that case, however, $\beta \leftarrow (\gamma \equiv \delta_i)$ occurs in a shorter semi-directed cycle $\gamma \to (\beta \equiv \delta_1) \Rightarrow \cdots \Rightarrow (\delta_i \equiv \gamma)$ in $\overline{\mathcal{G}^*}$, contradicting the minimality of $k$; hence, $\delta_i \neq \gamma$ for each $i$. By Lemma 3.10(a), the minimal-length semi-directed cycle has at least one line $\delta_j \stackrel{\text{w}}{-} \delta_{j+1} \in \overline{\mathcal{G}^*}$ ($1 \leq j \leq k-1$); consider the minimal such $j$.

Suppose first that $j \geq 2$, so $\delta_{j-1} \to \delta_j - \delta_{j+1}$ occurs as a subgraph of $\overline{\mathcal{G}^*}$. By Lemma 3.9(b), $\delta_{j-1} \Rightarrow \delta_{j+1} \in \overline{\mathcal{G}^*}$, hence $\alpha \to (\beta \equiv \delta_1) \Rightarrow \cdots \Rightarrow \delta_{j-1} \Rightarrow \delta_{j+1} \Rightarrow \cdots \Rightarrow (\delta_k \equiv \alpha)$ is a shorter semi-directed cycle in $\overline{\mathcal{G}^*}$ containing $\alpha \to \beta$, contradicting the minimality of $k$.

Suppose next that $j = 1$, so $(\beta \equiv \delta_1) - \delta_2 \in \overline{\mathcal{G}^*}$. By Lemma 3.9(b), $\alpha \Rightarrow \delta_2 \Leftarrow \gamma$ occurs as an induced subgraph in $\overline{\mathcal{G}^*}$. By Lemma 3.10(c), neither edge can be a line, hence $\alpha \to \delta_2 \leftarrow \gamma$ occurs as an immorality in $\overline{\mathcal{G}^*}$. But now $\alpha \to \delta_2$ occurs in a shorter semi-directed cycle $\alpha \to \delta_2 \Rightarrow \cdots \Rightarrow (\delta_k \equiv \alpha)$, contradicting the minimality of $k$.

(b) By Lemma 3.9(a), $\overline{\mathcal{G}^*}$ has no flags. Assume that $\alpha \to \beta - \gamma$ occurs as a flag in $(\overline{\mathcal{G}^*})^\circ$. Then necessarily $\alpha \to \beta \in \overline{\mathcal{G}^*}$, so by (a), $\alpha \to \beta \to \gamma$ occurs as a chordless 2-dipath in $\overline{\mathcal{G}^*}$, and $\beta \to \gamma$ occurs in some semi-directed cycle in $\overline{\mathcal{G}^*}$: $\beta \to (\gamma \equiv \delta_1) \Rightarrow \delta_2 \Rightarrow \cdots \Rightarrow (\delta_k \equiv \beta)$, with $k \geq 4$ by Lemma 3.10(c). Choose the flag $\alpha \to \beta - \gamma$ in $(\overline{\mathcal{G}^*})^\circ$ such that the associated semi-directed cycle has minimal length $k$. [Note that $\delta_i \neq \alpha$ for each $i$; otherwise $\alpha \to \beta$ would be included in a semi-directed cycle in $\overline{\mathcal{G}^*}$ and would thus be converted to $\alpha - \beta$ in $(\overline{\mathcal{G}^*})^\circ$, a contradiction to the original assumption.]

By Lemma 3.10(a), the semi-directed cycle must have at least one line $\delta_j - \delta_{j+1} \in \overline{\mathcal{G}^*}$ ($1 \leq j \leq k-1$); consider the minimal such $j$. If $j \geq 2$, then the minimality of $k$ is contradicted exactly as in (a), hence $j = 1$ and $\delta_1 - \delta_2 \in \overline{\mathcal{G}^*}$. By Lemma 3.9(b), $\beta \Rightarrow \delta_2 \in \overline{\mathcal{G}^*}$. If $\beta - \delta_2$, then $\beta \to (\gamma \equiv \delta_1) - \delta_2 - \beta$ is a semi-directed 3-cycle containing $\beta \to \gamma$ in $\overline{\mathcal{G}^*}$, contradicting Lemma 3.10(c); hence $\beta \to \delta_2 \in \overline{\mathcal{G}^*}$.

If $\alpha \cdot/\cdot \delta_2$ in $\overline{\mathcal{G}^*}$, then $\alpha \to \beta \to \delta_2$ occurs as a chordless 2-dipath in $\overline{\mathcal{G}^*}$ and $\beta \to \delta_2$ occurs in the semi-directed cycle $\beta \to \delta_2 \Rightarrow \cdots \Rightarrow (\delta_k \equiv \beta)$, contradicting the minimality of $k$. If $\alpha \cdots \delta_2$ in $\overline{\mathcal{G}^*}$, then $\alpha \to \delta_2$, since otherwise $\alpha \to \beta$ would occur in the semi-directed cycle $\alpha \to \beta \to \gamma - \delta_2 \Rightarrow \alpha$ in $\overline{\mathcal{G}^*}$, hence would be converted into $\alpha - \beta$ in $(\overline{\mathcal{G}^*})^\circ$. But now $\alpha \to \delta_2 - \gamma$ is a flag in $\overline{\mathcal{G}^*}$, contradicting Lemma 3.9(a).

(c) Suppose that $(\overline{\mathcal{G}^*})^\circ$ has a chordless undirected cycle which has no triplexes. By Lemma 3.9(a), $\overline{\mathcal{G}^*}$ has no such cycles, so has at least one arrow in this cycle. If $\overline{\mathcal{G}^*}$ has at least one opposing arrow in this cycle, it must have at least one triplex therein. If it has no opposing arrow, then it has at least one line [Lemma 3.10(a)], so again at least one triplex therein. But by (a) and (b), $\overline{\mathcal{G}^*}$ and $(\overline{\mathcal{G}^*})^\circ$ have the same triplexes, hence a contradiction. □

THEOREM 3.2. *$\mathcal{G}^*$ is adicyclic (i.e., is a chain graph) and $\mathcal{G}^* \in \mathcal{G}$.*

PROOF. By Lemma 3.8(b), to show that $\mathcal{G}^*$ is adicyclic, it suffices to show that $\overline{\mathcal{G}^*}$ is adicyclic, that is, $\overline{\mathcal{G}^*} = (\overline{\mathcal{G}^*})^\circ$. To simplify notation, set $K := \overline{\mathcal{G}^*}$. By Lemma 3.11, $K$ and $K^\circ$ have the same immoralities, have no flags and $K^\circ$ is a chain graph each of whose chain components $(K^\circ)_\eta$, $\eta \in \Xi(K^\circ)$, is chordal. [Note that $\eta \subseteq \Sigma \equiv \Sigma(\mathcal{G})$.]

Because $K \subseteq K^\circ$, it suffices to show that if $\alpha - \beta \in K^\circ$, then $\alpha - \beta \in K$. Let $\eta \in \Xi(K^\circ)$ be the unique chain component of $K^\circ$ such that $\alpha - \beta \in (K^\circ)_\eta$. Since it is chordal, $(K^\circ)_\eta$ admits two perfect directed versions, say $D_\eta$ and $D'_\eta$, such that $\alpha \to \beta \in D_\eta$ and $\alpha \leftarrow \beta \in D'_\eta$. (Apply Maximum Cardinality Search starting first at $\alpha$ and next at $\beta$.) Extend $D_\eta$ and $D'_\eta$ to directed graphs $D$ and $D'$, each having the same vertex set $\Sigma$ and the same skeleton as $K$ and $K^\circ$, by assigning perfect orientations (the same for $D$ and $D'$) to all other chain components $(K^\circ)_\chi$, $\chi \in \Xi(K^\circ)$, $\chi \neq \eta$. Thus, every arrow in $K^\circ$ also occurs in $D$ and $D'$, while $\alpha \to \beta \in D$ and $\alpha \leftarrow \beta \in D'$.

It is readily verified that $D$ and $D'$ are acyclic (since $K^\circ$ is adicyclic and perfect orientations of chordal graphs are acyclic) and have the same immoralities as $K$ and $K^\circ$ (since $K^\circ$ has no flags, perfect orientations of chordal graphs are moral, and every arrow in $K^\circ$ also occurs in $D$ and $D'$).

Now consider the "un-reduced" versions $H$ of $D$ and $H'$ of $D'$. That is, $H$ and $H'$ have the same skeleton as $\mathcal{G}^*$, while if $c \cdots d \in \mathcal{G}^*$, then:

A: $c - d \in H, H'$ iff $c \stackrel{s}{-} d \in \mathcal{G}^*$, that is, $H_\sigma = H'_\sigma = \mathcal{G}^*_\sigma\ \forall\sigma \in \Sigma$;
B: $c \to d \in H$ (resp., $H'$) iff $\sigma(c) \to \sigma(d) \in D$ (resp., $D'$).

In case B, either

B1: $\sigma(c) \to \sigma(d) \in K^\circ$, in which case $\sigma(c) \to \sigma(d) \in K$, or
B2a: $\sigma(c) - \sigma(d) \in K^\circ$ and $\sigma(c) \stackrel{w}{-} \sigma(d) \in K$, or
B2b: $\sigma(c) - \sigma(d) \in K^\circ$ and $\sigma(c) \leftrightarrow \sigma(d) \in K$ but the arrow $\leftrightarrow$ is contained in a semi-directed cycle in $K$.

Note that if $\tilde{a} \in \alpha$ and $\tilde{b} \in \beta$ are chosen such that $\tilde{a} \cdots \tilde{b} \in \mathcal{G}^*$, then $\tilde{a} \to \tilde{b} \in H$ and $\tilde{a} \leftarrow \tilde{b} \in H'$. If we can show that $H, H' \in \mathcal{G}$, then $\tilde{a} - \tilde{b} \in \mathcal{G}^*$, so $\alpha - \beta \in \overline{\mathcal{G}^*} \equiv K$, which would complete the proof of the adicyclicity of $\mathcal{G}^*$. Since $D$ and $D'$ are acyclic, $H$ and $H'$ are adicyclic by their construction. Therefore, it suffices to show that $H$ and $H'$ have the same immoralities and the same flags as $\mathcal{G}^*$, hence, the same triplexes.

(i) Suppose that $c \to d \leftarrow e$ occurs as an immorality in $\mathcal{G}^*$. By Lemma 3.5(a), this immorality also occurs in $(\mathcal{G}^*)^\circ$, so the arrow(s) $\sigma(c) \to \sigma(d)$ and $\sigma(d) \leftarrow \sigma(e)$ occur in $\overline{(\mathcal{G}^*)^\circ} = (\overline{\mathcal{G}^*})^\circ \equiv K^\circ$ [Lemma 3.8(b)] and therefore occur in $D$ and $D'$. [Note that it is possible that $\sigma(c) = \sigma(e)$.] Thus, $c \to d$ and $d \leftarrow e$ both occur in $H$ and $H'$ so, since $c \not\mathrel{\cdot/\cdot} e$ in $H$ and $H'$, $c \to d \leftarrow e$ occurs as an immorality in $H$ and $H'$.

(ii) Suppose that $c \stackrel{\mathrm{s}}{\to} d \stackrel{\mathrm{s}}{-} e$ occurs as a flag in $\mathcal{G}^*$. By Lemma 3.5(b), this flag also occurs in $(\mathcal{G}^*)^\circ$, so the arrow $\sigma(c) \to \sigma(d)$ occurs in $\overline{(\mathcal{G}^*)^\circ} = (\overline{\mathcal{G}^*})^\circ \equiv K^\circ$ and therefore in $D$ and $D'$. Thus, $c \to d - e$ occurs as a flag in $H$ and $H'$.

(iii) Suppose that $c \to d - e$ occurs as a flag in $H$ (or $H'$). Then $c \cdots d \stackrel{\mathrm{s}}{-} e$ occurs as an induced subgraph of $\mathcal{G}^*$ and $\sigma(c) \to \sigma(d) \in D$ (or $D'$). Consider the three possibilities B1, B2a, B2b for the edge $\sigma(c) \cdots \sigma(d)$:

B1: Here $\sigma(c) \to \sigma(d) \in K \equiv \overline{\mathcal{G}^*}$, so $c \to d \stackrel{\mathrm{s}}{-} e$ occurs as a flag in $\mathcal{G}^*$.

B2a: Here $\sigma(c) \stackrel{\mathrm{w}}{-} \sigma(d) \in K$ so $c \stackrel{\mathrm{w}}{-} d \in \mathcal{G}^*$, hence $c \stackrel{\mathrm{w}}{-} d \stackrel{\mathrm{s}}{-} e \in \mathcal{G}^*$, hence $c \stackrel{\mathrm{w}}{-} e \in \mathcal{G}^*$ [Lemma 3.2(b)]. But $c \not\mathrel{\cdot/\cdot} e$ in $\mathcal{G}^*$, a contradiction.

B2b: Here $\sigma(c) \leftrightarrow \sigma(d) \in K$ and the arrow occurs in a semi-directed cycle in $K \equiv \overline{\mathcal{G}^*}$, so by Lemma 3.8(a), $c \leftrightarrow d$ must occur in a semi-directed cycle in $\mathcal{G}^*$. Thus, either $c \to d \stackrel{\mathrm{s}}{-} e$ occurs as a flag in $\mathcal{G}^*$ but not in $(\mathcal{G}^*)^\circ$, contradicting Lemma 3.5(b), or $c \leftarrow d \stackrel{\mathrm{s}}{-} e$ occurs as an antiflag in $\mathcal{G}^*$ but not in $(\mathcal{G}^*)^\circ$, contradicting Lemma 3.5(c).

(iv) Suppose that $c \to d \leftarrow e$ occurs as an immorality in $H$ (or $H'$). Then $c \cdots d \cdots e$ occurs as an induced subgraph of $\mathcal{G}^*$, and $\sigma(c) \to \sigma(d) \in D$ (or $D'$), $\sigma(c) \to \sigma(d) \in D$ (or $D'$).

First, assume that $\sigma(c) = \sigma(e)$. Again consider the three possibilities B1, B2a, B2b for the edge $\sigma(c) \cdots \sigma(d)$:

B1: Here $\sigma(c) \to \sigma(d) \in K \equiv \overline{\mathcal{G}^*}$, so $c \to d \leftarrow e$ occurs as an immorality in $\mathcal{G}^*$.

B2a: Here $\sigma(c) \stackrel{\mathrm{w}}{-} \sigma(d) \in K$, so $c \stackrel{\mathrm{w}}{-} d \in \mathcal{G}^*$, $e \stackrel{\mathrm{w}}{-} d \in \mathcal{G}^*$, and $\mathcal{G}^*_{\sigma(c)} \equiv \mathcal{G}^*_{\sigma(e)}$ is complete [Lemma 3.6(d)], but $c \not\mathrel{\cdot/\cdot} e$ in $\mathcal{G}^*$, a contradiction.

B2b: Here $\sigma(c) \leftrightarrow \sigma(d) \in K$ and the arrow occurs in a semi-directed cycle in $K \equiv \overline{\mathcal{G}^*}$, so $\mathcal{G}^*_{\sigma(c)} \equiv \mathcal{G}^*_{\sigma(e)}$ is complete [Lemma 3.10(d)]. But $c \not\mathrel{\cdot/\cdot} e$ in $\mathcal{G}^*$, a contradiction.

Second, assume that $\sigma(c) \neq \sigma(e)$. As above, three possibilities (B1, B2a, B2b) exist for the edge $\sigma(c) \cdots \sigma(d)$ in $K^\circ$. Similarly, three possibilities (B1', B2a', B2b') exist for the edge $\sigma(e) \cdots \sigma(d)$ in $K^\circ$:

B1': $\sigma(e) \to \sigma(d) \in K^\circ$, in which case $\sigma(e) \to \sigma(d) \in K$, or

B2a': $\sigma(e) - \sigma(d) \in K^\circ$ and $\sigma(e) \stackrel{\mathrm{w}}{-} \sigma(d) \in K$, or

B2b': $\sigma(e) - \sigma(d) \in K^\circ$ and $\sigma(e) \leftrightarrow \sigma(d) \in K$, but the arrow $\leftrightarrow$ is contained in a semi-directed cycle in $K$.

Thus, we must consider the nine cases (B1, B1′), (B2a, B1′), ..., (B2b, B2b′).

(B1, B1′): Here $\sigma(c) \to \sigma(d) \in K$ and $\sigma(e) \to \sigma(d) \in K \equiv \overline{\mathcal{G}^*}$, so $c \to d \leftarrow e$ occurs as an immorality in $\mathcal{G}^*$.

(B1, B2a′): Here $\sigma(c) \to \sigma(d) \in K$ and $\sigma(e) \stackrel{\text{w}}{-} \sigma(d) \in K$, so $c \to d \stackrel{\text{w}}{-} e$ occurs as a flag in $\mathcal{G}^*$, contradicting Lemma 3.2(a). Similarly, the case (B2a, B1′) is also impossible.

(B1, B2b′): Here $\sigma(c) \to \sigma(d) \in K$ and $\sigma(e) \leftrightarrow \sigma(d) \in K$ with the arrow $\leftrightarrow$ contained in a semi-directed cycle in $K \equiv \overline{\mathcal{G}^*}$. By Lemma 3.8(a), $e \leftrightarrow d$ must occur in a semi-directed cycle in $\mathcal{G}^*$, so either $c \to d \leftarrow e$ occurs as an immorality in $\mathcal{G}^*$ but not in $(\mathcal{G}^*)^\circ$, contradicting Lemma 3.5(a), or $c \to d \to e$ occurs as a chordless 2-dipath in $\mathcal{G}^*$ but $c \to d - e$ occurs as a flag in $(\mathcal{G}^*)^\circ$, contradicting Lemma 3.5(b). Similarly, (B2b, B1′) is impossible.

In the remaining four cases (B2a, B2a′), (B2a, B2b′), (B2b, B2a′) and (B2b, B2b′), $\sigma(c) - \sigma(d) - \sigma(e)$ occurs as a subgraph in $K^\circ$. Because $\sigma(c) \to \sigma(d) \leftarrow \sigma(e)$ occurs as a subgraph in $D$ (or $D'$) and $K^\circ$ has the same immoralities as $D, D'$, necessarily $\sigma(c) \cdots \sigma(e) \in K \equiv \overline{\mathcal{G}^*}$ in these four cases. It cannot occur that $\sigma(c) \stackrel{\text{w}}{-} \sigma(e) \in K$ [otherwise $c \cdots e \in \mathcal{G}^*$ by Lemma 3.6(d)], so $\sigma(c) \leftrightarrow \sigma(e) \in K$. Therefore, $\exists c' \in \sigma(c), e' \in \sigma(e)$ such that $c' \leftrightarrow e' \in \mathcal{G}^*$. Since $c \cdot / \cdot e$ in $\mathcal{G}^*$, either $c \neq c'$ or $e \neq e'$ (or both).

(B2a, B2a′): Here $\sigma(c) \stackrel{\text{w}}{-} \sigma(d) \in K$ and $\sigma(e) \stackrel{\text{w}}{-} \sigma(d) \in K \equiv \overline{\mathcal{G}^*}$, so $c \stackrel{\text{w}}{-} d \in \mathcal{G}^*$, $e \stackrel{\text{w}}{-} d \in \mathcal{G}^*$, and both $\mathcal{G}^*_{\sigma(c)}$ and $\mathcal{G}^*_{\sigma(e)}$ are complete [Lemma 3.6(d)]. Thus, $c' \leftrightarrow e'$ occurs in a semi-directed cycle in $\mathcal{G}^*$.

(•) If we assume that $c \neq c'$, then $c \stackrel{\text{s}}{-} c' \leftrightarrow e'$ occurs as a subgraph of $\mathcal{G}^*$. In this case $c \leftrightarrow e' \in \mathcal{G}^*$, for otherwise the subgraph would occur as an antiflag or flag in $\mathcal{G}^*$ but not in $(\mathcal{G}^*)^\circ$, contradicting Lemma 3.5(c),(b). Therefore, $e' \neq e$ (since $c \cdot / \cdot e$ in $\mathcal{G}^*$), so $c \leftrightarrow e' \stackrel{\text{s}}{-} e$ occurs as a subgraph of $\mathcal{G}^*$. By the same reasoning, $c \leftrightarrow e \in \mathcal{G}^*$, thereby contradicting the fact that $c \cdot / \cdot e \in \mathcal{G}^*$. The same contradiction is obtained if we assume first that $e \neq e'$.

(B2a, B2b′): Here $\sigma(c) \stackrel{\text{w}}{-} \sigma(d) \in K$ and $\sigma(e) \leftrightarrow \sigma(d) \in K$ with the arrow $\leftrightarrow$ contained in some semi-directed cycle in $K$. Thus, $\mathcal{G}^*_{\sigma(c)}$ is complete by Lemma 3.6(d) and $\mathcal{G}^*_{\sigma(e)}$ is complete by Lemma 3.10(d), so $c \stackrel{\text{s}}{-} c' \in \mathcal{G}^*$ if $c \neq c'$ and $e \stackrel{\text{s}}{-} e' \in \mathcal{G}^*$ if $e \neq e'$. Because the triangle $\sigma(c) \stackrel{\text{w}}{-} \sigma(d) \leftrightarrow \sigma(e) \leftrightarrow \sigma(c)$ occurs in $K \equiv \overline{\mathcal{G}^*}$ and has exactly one line, by Lemma 3.10(b), it cannot occur as a semi-directed 3-cycle. This leaves two possible configurations for the triangle in $\overline{\mathcal{G}^*}$: $\sigma(c) \stackrel{\text{w}}{-} \sigma(d) \leftarrow \sigma(e) \to \sigma(c)$ and $\sigma(c) \stackrel{\text{w}}{-} \sigma(d) \to \sigma(e) \leftarrow \sigma(c)$, and in both cases the arrow between $\sigma(d)$ and $\sigma(e)$ is contained in some semi-directed cycle in $\overline{\mathcal{G}^*}$.

If the first configuration obtains, therefore, there exists a path $\sigma(d) \Rightarrow \delta_1 \Rightarrow \cdots \Rightarrow (\delta_k \equiv \sigma(e))$ in $\overline{\mathcal{G}^*}$ ($k \geq 2$). If each $\delta_i \neq \sigma(c)$, then $\sigma(e) \to \sigma(c)$ is contained in the semi-directed cycle $\sigma(e) \to \sigma(c) \stackrel{\text{w}}{-} \sigma(d) \Rightarrow \delta_1 \Rightarrow \cdots \Rightarrow (\delta_k \equiv \sigma(e))$ in $\overline{\mathcal{G}^*}$, while if $\delta_i = \sigma(c)$ for some (at most one) $i = 2, \ldots, k-1$ (so $k \geq 3$), then $\sigma(e) \to \sigma(c)$ is contained in the semi-directed cycle $\sigma(e) \to (\sigma(c) \equiv \delta_i) \Rightarrow \cdots \Rightarrow (\delta_k \equiv \sigma(e))$ in $\overline{\mathcal{G}^*}$. In both cases, therefore, $e' \to c'$ is contained in a semi-directed cycle in $\mathcal{G}^*$ by Lemma 3.8(a). If the second configuration obtains, a similar argument shows that $c' \to e'$ is contained in a semi-directed cycle in $\mathcal{G}^*$. Now a contradiction is reached exactly as in ($\bullet$). Similarly, (B2b, B2a$'$) is also impossible.

(B2b, B2b$'$): Here $\sigma(c) \leftrightarrow \sigma(d) \in K$ and $\sigma(e) \leftrightarrow \sigma(d) \in K$, with the arrows $\leftrightarrow$ both contained in semi-directed cycles in $K$. Since also $\sigma(c) \leftrightarrow \sigma(e) \in K$, $\sigma(c) \leftrightarrow \sigma(d) \leftrightarrow \sigma(e) \leftrightarrow \sigma(c)$ occurs as a triangle in $K \equiv \overline{\mathcal{G}^*}$. Of the eight possible orientations, two are cyclic but have no lines, hence are impossible by Lemma 3.4(b), leaving six acyclic possibilities: $\sigma(c) \to \sigma(d) \leftarrow \sigma(e) \to \sigma(c)$, $\sigma(c) \leftarrow \sigma(d) \leftarrow \sigma(e) \to \sigma(c)$, and so on. The argument in the preceding paragraph can be extended to show again that $c' \leftrightarrow e'$ occurs in a semi-directed cycle in $\mathcal{G}^*$, and again a contradiction is reached exactly as in ($\bullet$).

Thus, we have established that $H$ and $H'$ have the same triplexes as $\mathcal{G}^*$, so the proof of the adicyclicity of $\mathcal{G}^*$ is complete. Last, it follows now from Lemma 3.1 that $\mathcal{G}^* \in \mathcal{G}$. $\square$

REMARK 3.1. For subsequent use, we summarize the properties of the graphs $H$ and $H'$ constructed in the preceding proof: $H$ and $H'$ are adicyclic and have the same skeleton, immoralities and flags as $\mathcal{G}^*$, so $H, H' \in \mathcal{G}$. Each strong line in $\mathcal{G}^*$ occurs as a line in $H$ and $H'$, each arrow in $\mathcal{G}^*$ occurs as an arrow with the same orientation in $H$ and $H'$, and each weak line in $\mathcal{G}^*$ is converted to an arrow in $H$ and $H'$. Further, if a weak line $a \stackrel{\text{w}}{-} b \in \mathcal{G}^*$ is specified, $H$ and $H'$ can be chosen such that $a \to b \in H$ and $a \leftarrow b \in H'$.

**4. Local properties of AMP essential graphs.** Having established that $\mathcal{G}^*$ and $\overline{\mathcal{G}^*}$ are adicyclic, some of the results in Lemmas 3.2 and 3.9 can be sharpened and extended. These then yield information about the possible local configurations of strong/weak arrows/lines in $\mathcal{G}^* \equiv (V, E^*)$.

A chain component $\xi \in \Xi(\mathcal{G}^*)$ is *nontrivial* if $|\xi| \geq 2$, so that $\mathcal{G}^*_\xi$ contains at least one line. A nontrivial chain component $\xi$ is called *strong (weak)* if each line in $\mathcal{G}^*_\xi$ is strong (weak). For $v \in V$, denote the unique chain component of $\mathcal{G}^*$ that contains $v$ by $\xi(v) \equiv \xi_{\mathcal{G}^*}(v)$.

LEMMA 4.1. (a) *If $a \to b \stackrel{\text{w}}{-} c$ occurs as a subgraph in $\mathcal{G}^*$, then $a \to c \in \mathcal{G}^*$. Thus, if $\xi(b)$ is weak, then $a \to b' \, \forall \, b' \in \xi(b)$.*
(b) *If $\alpha \to \beta - \gamma$ occurs as a subgraph in $\overline{\mathcal{G}^*}$, then $\alpha \to \gamma \in \overline{\mathcal{G}^*}$.*

(c) *If $a \to b \overset{\mathrm{w}}{-} c$ occurs as a subgraph in $\mathcal{G}^*$, then $a \to b' \overset{\mathrm{w}}{-} c'$ and $a \to c'$ occur as subgraphs of $\mathcal{G}^*$ for all $b' \in \sigma(b)$ and $c' \in \sigma(c)$.*

(d) *The configuration $a \overset{\mathrm{s}}{-} b \overset{\mathrm{w}}{-} c$ cannot occur as an induced or noninduced subgraph in $\mathcal{G}^*$.*

PROOF. (a) and (b) By the adicyclicity of $\mathcal{G}^*$ and $\overline{\mathcal{G}^*}$, these results follow immediately from Lemmas 3.2(b) and 3.9(b) and the connectivity of $\xi(b)$.

(c) By Lemma 3.6(d), $a \to b \overset{\mathrm{w}}{-} c' \overset{\mathrm{w}}{-} b'$ occurs as a subgraph in $\mathcal{G}^*$. Therefore, $a \to c' \overset{\mathrm{w}}{-} b'$ occurs in $\mathcal{G}^*$ by (a), whence $a \to b' \in \mathcal{G}^*$ also by (a). [Note that if the arrow $a \to b$ is weak, the result concerning $a \to b' \overset{\mathrm{w}}{-} c'$ follows from Lemma 3.6(c), (d), and in turn implies that $a \to c' \in \mathcal{G}^*$ by (a).]

(d) Assume that $a \overset{\mathrm{s}}{-} b \overset{\mathrm{w}}{-} c$ occurs as a subgraph in $\mathcal{G}^*$ and let $\xi = \xi(a) \equiv \xi(b) \equiv \xi(c)$. By Lemma 3.6(d), $a' \overset{\mathrm{w}}{-} c' \in \mathcal{G}^* \ \forall a' \in \sigma(a), c' \in \sigma(c)$, and $\sigma(a) \equiv \sigma(b)$ and $\sigma(c)$ are complete subsets of $\xi$. Since every $\sigma \in \Sigma(\mathcal{G})$ such that $\sigma \subseteq \xi$ must be connected to $\sigma(a)$ by a path of weak and strong lines, it follows that every such $\sigma$ is complete. This implies that any chordless cycle $\mathcal{C}$ in $\mathcal{G}_\xi^*$ must contain at least one weak line $d \overset{\mathrm{w}}{-} e$, hence $\exists G \in \mathcal{G}$ such that either $d \to e \in G$ or $d \leftarrow e \in G$. But this and the adicyclicity of $G$ imply the existence of a triplex in $G_\mathcal{C}$, while $\mathcal{G}_\mathcal{C}^*$, being undirected, has no triplexes, a contradiction since $G$ and $\mathcal{G}^*$ have the same triplexes. Thus, $\mathcal{G}_\xi^*$ can have no chordless cycles, hence is chordal.

Thus, $\mathcal{G}_\xi^*$ admits a perfect directed version $F_\xi$ (apply MCS), so the edge $a \cdots b$ occurs as an arrow in $F_\xi$. Let $F$ be the graph obtained from $\mathcal{G}^*$ if we replace $\mathcal{G}_\xi^*$ by $F_\xi$, so $F \subseteq \mathcal{G}^*$. Clearly, $F$ has the same skeleton as $\mathcal{G}^*$ and is adicyclic: any semi-directed cycle in $F$ cannot be wholly contained in $\xi$ since $F_\xi$ is acyclic, hence such a cycle must contain at least one arrow from $\mathcal{G}^*$, so would correspond to a semi-directed cycle in $\mathcal{G}^*$ which, however, is adicyclic. Because $F$ differs from $\mathcal{G}^*$ only in that all lines in $\mathcal{G}_\xi^*$ become arrows in $F$, because no triplexes occur within $F_\xi$ or $\mathcal{G}_\xi^*$ (the orientation in $F_\xi$ is perfect, while $\mathcal{G}_\xi^*$ is undirected), and because neither $F$ nor $\mathcal{G}^*$ has a line $d - e$ with $e \in \xi$ and $f \in V \setminus \xi$ [since $\xi \in \Xi(\mathcal{G})$], the triplexes of $F$ and $\mathcal{G}^*$ can differ only if, for some $e, f \in \xi$, $d \to e - f$ occurs as a flag in $\mathcal{G}^*$ but $d \to e \to f \in F$. The former is impossible if $e - f$ is weak in $\mathcal{G}^*$ [Lemma 3.2(a)], so assume $e \overset{\mathrm{s}}{-} f \in \mathcal{G}^*$. Since $\sigma(e) \subsetneq \xi$, $\exists g \in (\xi \setminus \sigma(e))$ such that $e \overset{\mathrm{w}}{-} g \in \mathcal{G}^*$ [Lemma 3.6(d)], so, since $f \in \sigma(e)$, $d \to f \in \mathcal{G}^*$ [set $a = d$, $b = e$, $c = g$ and $b' = f$ in part (c)]. Thus, $d \to e - f$ cannot occur as a flag in $\mathcal{G}^*$, so $F$ and $\mathcal{G}^*$ have the same triplexes, hence $F \in \mathcal{G}$. Since $a \cdots b$ occurs as an arrow in $F$, $a - b$ cannot be a strong line in $\mathcal{G}^*$, contrary to assumption. □

PROPOSITION 4.1. (a) *Every nontrivial chain component of $\mathcal{G}^*$ is either strong or weak.*

(b) *A nontrivial chain component $\xi \in \Xi(\mathcal{G}^*)$ is strong iff $\mathcal{G}^*_{\bar{\xi}}$ contains either a chordless undirected cycle or a flag, where $\bar{\xi} = \mathrm{cl}_{\mathcal{G}^*}(\xi) \equiv \xi \dot{\cup} \mathrm{pa}_{\mathcal{G}^*}(\xi)$. In fact, $\mathcal{G}^*_{\bar{\xi}}$ must contain either a chordless undirected cycle or a biflag.*

PROOF. Lemma 4.1(d) yields (a); Lemmas 3.2(a) and 3.3 yield "if" in (b).

(b) "only if": Assume that $\xi$ is strong but $\mathcal{G}^*_{\bar{\xi}}$ contains no chordless undirected cycles or biflags. Suppose first that $\mathcal{G}^*_{\bar{\xi}}$ has no flags. Because $\mathcal{G}^*_\xi$ is chordal, we can define $F_\xi$ and $F$ as in the proof of Lemma 4.1(d). Then as above, $F$ is adicyclic and has the same skeleton and triplexes as $\mathcal{G}^*$, hence $F \in \mathcal{G}$. Because all edges in $F_\xi$ are arrows, however, the corresponding edges in $\mathcal{G}^*_\xi$ cannot be strong lines, a contradiction.

Suppose next that $\mathcal{G}^*_{\bar{\xi}}$ has at least one flag, say, $a^* \to b^* - c^*$ with $b^*, c^* \in \xi$. Let $F$ be the graph constructed from $\mathcal{G}^*$ by converting each such flag into an immorality $a^* \to b^* \leftarrow c^*$. This process is unambiguous since $b^* - c^*$ cannot occur in a 2-biflag in $\mathcal{G}^*$. Since $F \subseteq \mathcal{G}^*$ and they have the same skeleton, each chain component $\rho$ of $F$ satisfies either $\rho \subseteq \xi$ or $\rho \cap \xi = \varnothing$. Define

$$\Xi_\xi(F) := \{\rho \in \Xi(F) | \rho \subseteq \xi\},$$

so $\xi = \dot{\cup} \{\rho | \rho \in \Xi_\xi(F)\}$. We shall establish four properties of $F$:

(1) $F$ is adicyclic. If $F$ were not adicyclic, it would contain a semi-directed cycle $c_0 \to c_1 \Rightarrow \cdots \Rightarrow c_k \equiv c_0$ $(k \geq 3)$. Because $F \subseteq \mathcal{G}^*$ and $\mathcal{G}^*$ is adicyclic, this cycle must be completely undirected in $\mathcal{G}^*$, in particular, $c_0 - c_1 \in \mathcal{G}^*$, so there must exist a vertex $a \neq c_0, c_1, c_2, c_{k-1}$ such that $c_0 - c_1 \leftarrow a$ occurs as a flag in $\mathcal{G}^*$. (If $k \geq 5$, then also $a \neq c_i$ for $i = 3, \ldots, k-2$, for otherwise $c_1 \Rightarrow c_2 \cdots \Rightarrow c_i \to c_1$ would be a semi-directed cycle in $\mathcal{G}^*$.) Further, $a \to c_1 - c_2$ cannot occur as a flag in $\mathcal{G}^*$ (otherwise $c_1 \leftarrow c_2 \in F$), so $a \to c_2 \in \mathcal{G}^*$. Similarly, $a \to c_i \in \mathcal{G}^*$ for $i = 3, \ldots, k$, which contradicts the nonadjacency of $a$ and $c_k \equiv c_0$ in the flag $c_0 - c_1 \leftarrow a$.

(2) Each triplex in $\mathcal{G}^*$ corresponds to a triplex in $F$. This is immediate by the construction of $F$ from $\mathcal{G}^*$.

(3) Each immorality in $F$ corresponds to a triplex in $\mathcal{G}^*$ (but not every flag in $F$ need correspond to a triplex in $\mathcal{G}^*$). This holds since $\mathcal{G}^*_{\bar{\xi}}$ contains no 3-biflag chains.

(4) The *moralized* graph (cf. [4], page 40) $L \equiv L(\rho) := (F_{\mathrm{cl}_F(\rho)})^{\mathrm{m}}$ is chordal for each $\rho \in \Xi_\xi(F)$. If not, then $L$ would contain a chordless $k$-cycle $c_0 - c_1 - c_2 - \cdots - c_k \equiv c_0$ $(k \geq 4)$, where each $c_i \in \mathrm{cl}_F(\rho) \equiv \rho \dot{\cup} \mathrm{pa}_F(\rho)$. Since $L_\rho = F_\rho = \mathcal{G}^*_\rho$ and $\mathcal{G}^*_\xi$ is chordal, this cycle cannot lie entirely within $\rho$, so must intersect $\mathrm{pa}_F(\rho)$. Since $\mathrm{pa}_F(\rho)$ is complete in $L$, either (a) exactly one vertex of the cycle lies in $\mathrm{pa}_F(\rho)$, say, $c_1$, or (b) exactly two vertices of the cycle, necessarily consecutive, lie in $\mathrm{pa}_F(\rho)$, say, $c_0 \equiv c_k$ and $c_1$.

In case (a), $c_2, c_k \in \rho$, while $c_1 \in \text{pa}_F(\rho)$, so $c_1 \to c_2 \in F$ and $c_1 \to c_k \in F$ by (1) and the connectedness of $F_\rho$, hence $c_1 \Rightarrow c_2 \in \mathcal{G}^*$ and $c_1 \Rightarrow c_k \in \mathcal{G}^*$. Because $c_2 - c_3 \in F_\rho = \mathcal{G}^*_\rho$ and $c_{k-1} - c_k \in F_\rho = \mathcal{G}^*_\rho$, $c_1 - c_2 \in \mathcal{G}^*$ and $c_1 - c_k \in \mathcal{G}^*$ by the definition of $F$, so the cycle lies in $\mathcal{G}^*_\xi$ and is chordless, contradicting the chordality of $\mathcal{G}^*_\xi$.

In case (b), a similar argument shows that $c_1 \to c_2 \in F$ and $c_k \to c_{k-1} \in F$, while $c_1 - c_2 \in \mathcal{G}^*$ and $c_k - c_{k-1} \in \mathcal{G}^*$. Therefore, all vertices $c_1, \ldots, c_k$ lie in $\xi$, so $c_1$ and $c_k$ cannot be adjacent in $\mathcal{G}^*$ since $\mathcal{G}^*_\xi$ is assumed to be chordal. (Instead, the line $c_1 - c_k \equiv c_0$ in $L$ must have been added by moralization.) Furthermore, by the definition of $F$, there must exist vertices $a, b$ (possibly $a = b$) such that $c_1 - c_2 \leftarrow a$ and $b \to c_{k-1} - c_k$ occur as flags in $\mathcal{G}^*$. Thus, $a \neq c_1, c_2, c_3$ and, since the cycle is chordless in $L$, $a \neq c_4, \ldots, c_k$. Similarly, $b \neq c_1, \ldots, c_k$. Because $c_i - c_{i+1} \in F \subseteq \mathcal{G}^*$ for $i = 2, \ldots, k-2$ and $\mathcal{G}^*_\xi$ contains no biflag chains, necessarily $a \to c_i \in \mathcal{G}^*$ for $i = 3, \ldots, k-1$ and $b \to c_i \in \mathcal{G}^*$ for $i = k-2, \ldots, 2$. Thus, at least one of $[a; c_1, \ldots, c_k]$, $[b; c_1, \ldots, c_k]$ or $[a, b; c_1, \ldots, c_k]$ must occur as a $k$-biflag in $\mathcal{G}^*$, contradicting the assumption that $\mathcal{G}^*_\xi$ does not contain a biflag. Thus, (4) holds.

Now construct a graph $F' \subseteq F$ with the same skeleton as $F$, as follows. For each $\rho \in \Xi_\xi(F)$, assign a perfect orientation to the edges of the chordal undirected graph $L(\rho) \equiv (F_{\text{cl}_F(\rho)})^m$ according to MCS starting at an arbitrary vertex in $\text{pa}_F(\rho)$ (if any), obtaining a perfect digraph $D(\rho)$ with vertex set $\text{cl}_F(\rho)$. Because $\text{pa}_F(\rho)$ is complete in $L(\rho)$, MCS can be chosen to visit each vertex in $\text{pa}_F(\rho)$ before visiting any vertex in $\rho$ itself. Thus, for any edge between a vertex in $\text{pa}_F(\rho)$ and a vertex in $\rho$, the orientations of this edge in $D(\rho)$ and in $F$ (and in $F'$ defined below) are identical.

Let $F' \subseteq F$ be the graph obtained from $F$ by orienting each undirected edge in $F_\rho$ ($= \mathcal{G}^*_\rho$) according to its orientation in $D(\rho)$, for each $\rho \in \Xi_\xi(F)$. Note that, for any pair $d, e \in V$ with either $d \notin \xi$ or $e \notin \xi$ (or both), the edge $d \cdots e$ (if any) must be identical in $\mathcal{G}^*$, $F$ and $F'$. Also, *every* line in $\mathcal{G}^*_\xi$ is converted to an arrow in $F'$. We shall show that $F'$ is adicyclic and has the same triplexes as $\mathcal{G}^*$, hence $F' \in \mathcal{G}$. This contradicts the fact that every line in $\mathcal{G}^*_\xi$ (in particular, $b^* - c^*$) is strong, which will complete the proof.

If $F'$ were not adicyclic, it would have a semi-directed cycle which cannot lie entirely within any one $\rho \in \Xi_\xi(F)$ since each $D(\rho)$ is perfect. Nor can this cycle lie entirely outside $\xi$, since there $F'$ coincides with $F$ which is adicyclic. Thus, this cycle must contain two consecutive vertices $c_1 \cdots c_2$ such that either $c_1 \in \rho_1$ and $c_2 \in \rho_2$ for distinct $\rho_1, \rho_2 \in \Xi_\xi(F)$, or else $c_1 \in \xi$ and $c_2 \notin \xi$. In either case the edge $c_1 \cdots c_2$ must occur as an arrow in $F$, so the cycle is also semi-directed in $F$ (since $F \supseteq F'$), contradicting the adicyclicity of $F$.

Since $\mathcal{G}^* \supseteq F'$, each immorality in $\mathcal{G}^*$ remains an immorality in $F'$, while, by the construction of $F$ and $F'$, each flag $d \to e - f$ in $\mathcal{G}^*$ either remains a flag in $F'$ (if $e, f \notin \xi$) or else is converted to an immorality in $F'$ (if $e, f \in \xi$). Also, by this construction, if $d \to e - f$ is a flag in $F'$, then necessarily $e, f \notin \xi$, so this flag also occurs in $\mathcal{G}^*$.

Last, suppose that $d \to e \leftarrow f$ is an immorality in $F'$ that does not occur as a triplex in $\mathcal{G}^*$. Then $d - e - f$ must occur as an induced subgraph in $\mathcal{G}^*_\xi$, but $d, e, f$ cannot all lie in the same $\rho \in \Xi_\xi(F)$, since $D(\rho)$ is perfect, hence has no immoralities. Without loss of generality, assume that $d \in \rho_1$ and $e \in \rho_2$, where $\rho_1 \neq \rho_2$. Thus, the edge $d \cdots e$ must be an arrow in $F$, hence $d \to e \in F$ since $F \supseteq F'$, and $e \Leftarrow f \in F$ for the same reason. If $e \leftarrow f \in F$, then $d \to e \leftarrow f$ would occur as an immorality in $F$, hence by (3) would correspond to a triplex in $\mathcal{G}^*$, contrary to assumption, so necessarily $e - f \in F$ and $f \in \rho_2$. Thus, the edge $d \cdots f$ cannot be added in the moralization process for $L(\rho_2)$, hence cannot occur in $D(\rho_2)$, so $d \to e \leftarrow f$ is an immorality in $D(\rho_2)$, contradicting its perfectness. Thus, $F'$ and $\mathcal{G}^*$ have the same triplexes. □

We note that not every line in a strong chain component need be contained in a biflag or chordless cycle. Examples appear in [5, 6].

Let $D_0 \equiv (V, E_0)$ be an ADG and let $\mathcal{D}$ denote its ADG Markov equivalence class, the set of all ADGs $D \equiv (V, E_D)$ that are Markov equivalent to $D_0$. Andersson, Madigan and Perlman [3] defined the *ADG essential graph*

(4.1) $$\mathcal{D}^* \equiv \bigcup \{D | D \in \mathcal{D}\} := (V, \bigcup (E_D \mid D \in \mathcal{D}))$$

determined by $\mathcal{D}$ and showed that it uniquely represents $\mathcal{D}$. For ADGs, in fact, the ADG and AMP definitions of essential graph are identical.

PROPOSITION 4.2. *$\mathcal{G}$ contains some ADG $D_0$ iff $\mathcal{G}^*$ has no biflags and all its chain components are chordal. In this case $\mathcal{G}^* = \mathcal{D}^*$.*

PROOF. Since $\mathcal{G}^* \in \mathcal{G}$, the first statement follows from Proposition 3 of [4]. By Proposition 4.1(b), $\mathcal{G}^*$ has no strong chain components, hence no strong lines. Clearly, each arrow in $\mathcal{G}^*$, strong or weak, must occur with the same orientation in every $D \in \mathcal{D}$, so must occur in $\mathcal{D}^*$. If a weak line $a \overset{\text{w}}{-} b$ occurs in $\mathcal{G}^*$, then the graphs $H, H'$ in Remark 3.1 are ADGs and belong to $\mathcal{D}$, so the line $a - b$ must occur in $\mathcal{D}^*$. Thus, $\mathcal{G}^* = \mathcal{D}^*$. □

We now turn to a characterization of the arrows in $\mathcal{G}^*$.

DEFINITION 4.1. Let $G \equiv (V, E)$ be a chain graph. An arrow $a \to b \in G$ is *protected* in $G$ if it occurs in at least one of the seven configurations shown in Figure 5 as an induced subgraph of $G$.

DEFINITION 4.2. Let $G \equiv (V, E)$ be a chain graph. An arrow $a \to b \in G$ is *irreversible* in $G$ if replacing $a \to b$ by $a \leftarrow b$ creates or destroys a triplex or creates a semi-directed cycle.

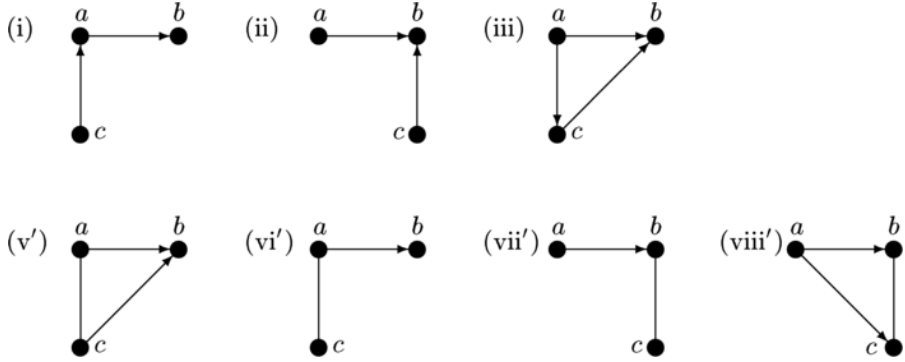

FIG. 5. *The seven protected configurations for an arrow $a \to b \in G$.*

Determination of the irreversibility of an arrow in $G$ apparently requires global knowledge of $G$, since semi-directed cycles may be of arbitrary length. In fact, however, only local knowledge of $G$ is required:

LEMMA 4.2. *Let $G \equiv (V, E)$ be a chain graph. An arrow is irreversible in $G$ if and only if it is protected in $G$.*

PROOF. Clearly, a protected arrow is irreversible. If $a \to b$ is irreversible in $G$ by virtue of $a \leftarrow b$ creating (resp., destroying) a triplex, then $a \to b$ must occur in configuration (i) or (vi′) [resp., (ii) or (vii′)] as an induced subgraph of $G$. If $a \to b$ is irreversible in $G$ by virtue of $a \leftarrow b$ creating a semi-directed cycle, then $a \to (b \equiv d_1) \Leftarrow d_2 \Leftarrow \cdots \Leftarrow (d_k \equiv a)$ $(k \geq 3)$ occurs as a subgraph in $G$. Since $G$ is adicyclic, at least one $\Leftarrow$ must be $\leftarrow$.

Suppose first that $b \leftarrow d_2 \in G$. If $a \cdot\!/\!\cdot d_2$, then $a \to b$ occurs in configuration (ii) with $c = d_2$. If $a \cdots d_2$ in $G$, then either $a \to d_2 \in G$, $a - d_2 \in G$ or $a \leftarrow d_2 \in G$. In the first (resp., second) case $a \to b$ occurs in configuration (iii) [resp., (v′)] with $c = d_2$. The third case cannot occur, for otherwise $k \geq 4$ and $a \leftarrow d_2 \Leftarrow \cdots \Leftarrow (d_k \equiv a)$ would be a semi-directed cycle in $G$. Suppose next that $b - d_2 \in G$. Then a similar argument shows that $a \to b$ must occur in configuration (vii′) or (viii′) with $c = d_2$. □

By the definition of $\mathcal{G}^*$, each arrow $a \to b \in \mathcal{G}^*$ must be irreversible in $\mathcal{G}^*$, hence, since $\mathcal{G}^*$ is a chain graph, protected in $\mathcal{G}^*$. In fact, each arrow must be *well protected* in $\mathcal{G}^*$.

DEFINITION 4.3. Let $\mathcal{G}^* \equiv (V, E^*)$ be an AMP essential graph. An arrow $a \to b \in \mathcal{G}^*$ is *well protected* in $\mathcal{G}^*$ if it occurs in at least one of the eight configurations shown in Figure 6 as an induced subgraph of $\mathcal{G}^*$.

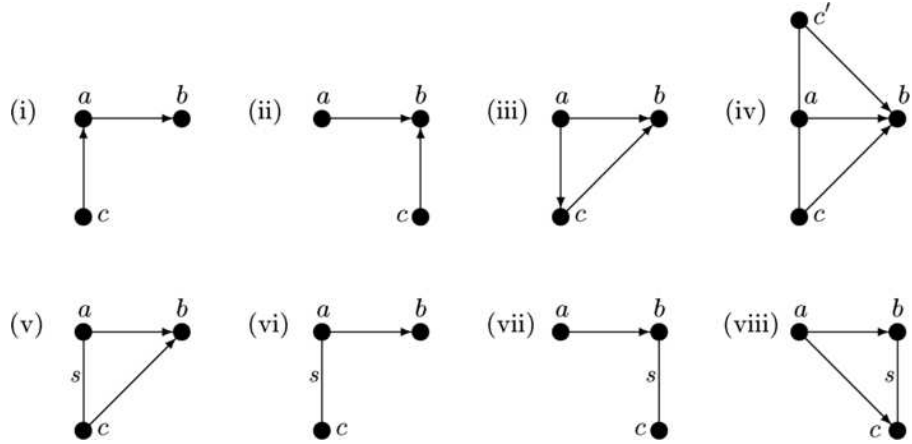

FIG. 6. *The eight well-protected configurations for an arrow $a \to b \in \mathcal{G}^*$*

PROPOSITION 4.3. *Each arrow in an AMP essential graph $\mathcal{G}^*$ is well protected in $\mathcal{G}^*$.*

PROOF. If $\xi(a)$ is strong, hence, nontrivial, $\exists c \in \xi(a)$ such that $c \stackrel{\text{s}}{-} a \to b$ is a subgraph of $\mathcal{G}^*$. Here $a \to b$ must occur in configuration (v) or (vi) in $\mathcal{G}^*$. If $\xi(b)$ is strong, hence nontrivial, then $\exists c \in \xi(a)$ such that $a \to b \stackrel{\text{s}}{-} c$ is a subgraph of $\mathcal{G}^*$. Here $a \to b$ must occur in configuration (vii) or (viii) in $\mathcal{G}^*$.

Last, assume that both $\xi(a)$ and $\xi(b)$ are weak (or trivial). In particular, $a \to b$ cannot occur in configuration (v), (vi), (vii) or (viii) in $\mathcal{G}^*$. We shall assume also that it does not occur in configuration (i), (ii), (iii) or (iv) and obtain a contradiction.

We begin by showing that $\mathcal{G}^*_\theta$ is complete, where

(4.2) $$\theta := \{c \in \xi(a) | c \to b \in \mathcal{G}^*\}.$$

Note that $a \in \theta$. If $c \in \theta \setminus \{a\}$, then $c \to b \in \mathcal{G}^*$. Thus, $c \cdots a$ (necessarily $c - a \in \mathcal{G}^*$), for otherwise $a \to b \leftarrow c$ would occur as an immorality in $\mathcal{G}^*$, contradicting the assumed nonoccurrence of $a \to b$ in configuration (ii) in $\mathcal{G}^*$. If $c, c' \in \theta \setminus \{a\}$ then $c \to b \in \mathcal{G}^*$ and $c' \to b \in \mathcal{G}^*$, while $c - a \in \mathcal{G}^*$ and $c' - a \in \mathcal{G}^*$ by the preceding argument, hence $c - c' \in \mathcal{G}^*$ by the nonoccurrence of (iv) in $\mathcal{G}^*$. Thus, $\theta$ is complete.

By Lemma 3.3, $\mathcal{G}^*_{\xi(b)}$ and $\mathcal{G}^*_{\xi(a)}$ are chordal UGs. Construct $G$ from $\mathcal{G}^*$ by assigning perfect orientations to the lines (if any) in these two UGs as follows. First, let $q = |\theta| \geq 1$, let $(c_1, \ldots, c_{q-1})$ be an arbitrary numbering of $\theta \setminus \{a\}$ and let $c_q = a$. Apply MCS to $\mathcal{G}^*_{\xi(a)}$ starting at $c_1$. The completeness of $\mathcal{G}^*_\theta$ ensures that MCS can reproduce the initial sequence $(c_1, \ldots, c_q \equiv a)$. The resulting perfect orientation in $G$ of the lines of $\mathcal{G}^*_{\xi(a)}$ satisfies the following two conditions:

($\alpha$) any line $a - c \in \mathcal{G}^*_{\xi(a)}$ with $c \in \theta \setminus \{a\}$ becomes $a \leftarrow c$ in $G$;

($\beta$) any line $a - d \in \mathcal{G}^*_{\xi(a)}$ with $d \in \xi(a) \setminus \theta$ becomes $a \rightarrow d$ in $G$. Next, orient the lines in $\mathcal{G}^*_{\xi(b)}$ according to MCS started at $b$; in particular,

($\gamma$) any line $b - c \in \mathcal{G}^*_{\xi(b)}$ becomes $b \rightarrow c$ in $G$.

Now construct $G'$ from $G$ by replacing $a \rightarrow b$ by $a \leftarrow b$, so $\mathcal{G}^*$, $G$ and $G'$ have the same skeleton. We shall show that $G$ and $G'$ are adicyclic and have the same flags and immoralities as $\mathcal{G}^*$, hence $G, G' \in \mathcal{G}$. Thus, $a - b \in \mathcal{G}^*$, the desired contradiction.

Because $G_{\xi(a)}$ and $G_{\xi(b)}$ are acyclic digraphs, any semi-directed cycle in $G$ cannot be entirely contained in $G_{\xi(a)}$ or in $G_{\xi(b)}$, hence must include at least one arrow of $\mathcal{G}^*$. Since $G \subseteq \mathcal{G}^*$, this cycle must correspond to a semi-directed cycle in $\mathcal{G}^*$, contradicting the adicyclicity of $\mathcal{G}^*$. Thus, $G$ is adicyclic.

If $d \rightarrow e \leftarrow f$ is an immorality in $\mathcal{G}^*$, then it also occurs as an immorality in $G$, since both graphs have the same skeleton and $G \subseteq \mathcal{G}^*$. If $d \rightarrow e - f$ is a flag in $\mathcal{G}^*$, then $d \rightarrow e \in G$ and $e \stackrel{s}{-} f \in \mathcal{G}^*$ [Lemma 3.2(a)] so $e - f$ remains a line in $G$, hence $d \rightarrow e - f$ is a flag in $G$. If $d \rightarrow e - f$ is a flag in $G$ but not in $\mathcal{G}^* \supseteq G$, then $d - e - f$ occurs as an induced subgraph of $\mathcal{G}^*$ and necessarily $d - e \in \mathcal{G}^*_{\xi(a)}$ or $\mathcal{G}^*_{\xi(b)}$, hence also $e - f \in \mathcal{G}^*_{\xi(a)}$ or $\mathcal{G}^*_{\xi(b)}$, so $e \cdots f$ must occur as an arrow in $G$, a contradiction. Finally, if $d \rightarrow e \leftarrow f$ is an immorality in $G$, but not in $\mathcal{G}^*$, then at least one of these arrows, say, $d \rightarrow e$, occurs as a line $d - e$ in $\mathcal{G}^*$, hence $d - e \in \mathcal{G}^*_{\xi(a)}$ or $\mathcal{G}^*_{\xi(b)}$, so $d \stackrel{w}{-} e \in \mathcal{G}^*$. Since $d \stackrel{w}{-} e \leftarrow f$ cannot occur as a flag in $\mathcal{G}^*$, necessarily $e - f \in \mathcal{G}^*$, so $d - e - f$ occurs as an induced subgraph of $\mathcal{G}^*_{\xi(b)}$ or $\mathcal{G}^*_{\xi(b)}$. But this subgraph cannot become an immorality $d \rightarrow e \leftarrow f$ in $G$ because the orientations of $G_{\xi(b)}$ and $G_{\xi(b)}$ are perfect, a contradiction. Thus, $\mathcal{G}^*$ and $G$ have the same flags and immoralities.

If $G'$ were not adicyclic, then since $G$ is adicyclic, $G'$ must contain a semi-directed cycle with $a \leftarrow b$, so this cycle must have the form $a \leftarrow (b \equiv c_1) \Leftarrow c_2 \Leftarrow \cdots \Leftarrow (c_k \equiv a)$ $(k \geq 3)$ in $G'$. Therefore, $G$ must contain the subgraph $a \rightarrow (b \equiv c_1) \Leftarrow c_2 \Leftarrow \cdots \Leftarrow (c_k \equiv a)$ and, since $G \subseteq \mathcal{G}^*$, a subgraph of this form also occurs in $\mathcal{G}^*$. Consider the edge $b \Leftarrow c_2$ in $\mathcal{G}^*$. If $b - c_2 \in \mathcal{G}^*$, then $b - c_2 \in \mathcal{G}^*_{\xi(b)}$, hence $b \rightarrow c_2 \in G$ by ($\gamma$), contradicting the occurrence of $b \Leftarrow c_2$ in $G$. Next, if $b \leftarrow c_2 \in \mathcal{G}^*$, then $a \Leftarrow c_2$ in $\mathcal{G}^*$ by the assumed nonoccurrence of (ii) and (iii) in $\mathcal{G}^*$. If this edge is $a \leftarrow c_2$ in $\mathcal{G}^*$, then $a \leftarrow c_2 \in G$, so $a \leftarrow c_2 \in G'$, which implies that $k \geq 4$ and that $a \leftarrow c_2 \Leftarrow \cdots \Leftarrow (c_k \equiv a)$ is a semi-directed cycle in $G$, a contradiction. If this edge is $a - c_2$ in $\mathcal{G}^*$, then $c_2 \in \theta \setminus \{a\}$, so $a \leftarrow c_2 \in G$ by ($\alpha$), leading to the same contradiction. Thus, $G'$ is adicyclic.

Suppose that an immorality $d \rightarrow e \leftarrow f$ occurs in $G$ but not in $G'$. Then one of these two arrows, say, $d \rightarrow e$, must be $a \rightarrow b$, so $a \rightarrow b \leftarrow f$ occurs as an immorality in $G \subseteq \mathcal{G}^*$. Thus, either $a \rightarrow b \leftarrow f$ occurs as an immorality in $\mathcal{G}^*$, contradicting the nonoccurrence of (ii), or $a \to b - f$ occurs as a subgraph in $\mathcal{G}^*$ so $b \to f \in G$ by $(\gamma)$, a contradiction. If $d \to e - f$ occurs as a flag in $G$ but not in $G'$, then the arrow $d \to e$ must be the arrow $a \to b$ in $G$, hence $a \to b - f$ occurs as a subgraph in $G$ and therefore in $\mathcal{G}^*$, so $b \to f \in G$ by $(\gamma)$, again a contradiction. If $d - e \leftarrow f$ occurs as a flag in $G'$ but not in $G$, then the arrow $e \leftarrow f$ must be the arrow $a \leftarrow b$ in $G'$, so $d - a \to b$ is a subgraph of $G \subseteq \mathcal{G}^*$, hence $d \in \xi(a)$. But then the edge $d \cdots a$ must be an arrow in $G$, a contradiction. Finally, if $d \to e \leftarrow f$ occurs as an immorality in $G'$ but not in $G$, then one arrow, say, $e \leftarrow f$, must be $a \leftarrow b$ in $G'$, so $d \to a \to b$ occurs as an induced subgraph in $G \subseteq \mathcal{G}^*$. Thus, either $d \to a \to b$ occurs as an induced subgraph in $\mathcal{G}^*$, contradicting the nonoccurrence of (i), or $d \overset{\text{w}}{-} a \to b$ occurs as an induced subgraph in $\mathcal{G}^*$, so $d \in \xi(a) \setminus \theta$, hence $d \leftarrow a \in G$ by $(\beta)$, a contradiction. Thus, $G$ and $G'$ have the same flags and immoralities. $\square$

The examples in Figures 3 and 4 show that Proposition 4.3 is rendered invalid if any of the eight configurations (i), (ii), (iii), (iv), (v), (vi), (vii) or (viii) is excluded from the definition of well protection in Definition 4.3. However, weak arrows and strong arrows of $\mathcal{G}^*$ are well protected in more stringent fashions, namely, configurations (iv) and (vii) cannot occur for the well protection of weak arrows, while (vi) is not required for the well protection of strong arrows—see [5], Propositions 4.4 and 4.5.

By Lemma 3.6(b), both arrows in configurations (v) and (viii) must be strong or both weak. The example in Figure 4(4) (resp., Figure 7) shows that if the strong line in (v) [resp., (viii)] is replaced by a weak line, then it can occur in $\mathcal{G}^*$ that one arrow is strong and the other weak.

We conclude this section with an easy characterization of those *directed* graphs that can occur as AMP essential graphs.

THEOREM 4.1. *A directed graph $G$ is an AMP essential graph iff it is acyclic and each arrow is protected in $G$, that is, occurs in one of the configurations* (i), (ii), *or* (iii) *as an induced subgraph of $G$. Thus, the class*

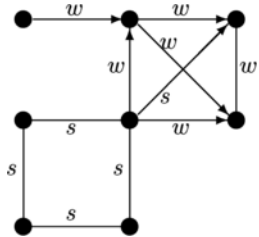

FIG. 7. *An AMP essential graph $\mathcal{G}^*$.*

*of directed AMP essential graphs coincides with the class of directed ADG essential graphs.*

PROOF. "Only if" follows immediately from Theorem 3.2 and Proposition 4.3. "If" and the final statement follow from Corollary 4.2 of [3] and Proposition 4.2. □

**5. Characterization of AMP essential graphs.** Because the well-protected configurations (v), (vi), (vii) and (viii) in Definition 4.3 involve strong lines, the local characterization of the arrows of an AMP essential graph $\mathcal{G}^*$ given in Proposition 4.3 in turn requires a characterization of the strong lines of $\mathcal{G}^*$. By Lemma 3.3, any line in a chordless undirected cycle or biflag must be strong, but $\mathcal{G}^*$ may contain other strong lines as well (see [5, 6]). Because [Proposition 4.1(b)] the strong chain components of $\mathcal{G}^*$ are distinguished from the weak chain components by the presence of at least one chordless cycle or biflag and the latter involves at least one arrow, the strong lines cannot be determined in an intrinsic way. A characterization (necessarily nonintrinsic) of the strong chain components of $\mathcal{G}^*$ is contained in the characterization of AMP essential graphs given in Theorem 5.1.

Some additional terminology is required. Let $G \equiv (V, E)$ be an adicyclic graph. For any nonempty subset $\alpha \subseteq V$, a vertex $v \notin \alpha$ is called a *covering neighbor* (*covering parent*) of $\alpha$ if $v$ is a neighbor (parent) of each $a \in \alpha$, while $v \notin \alpha$ is called a *noncovering neighbor* (*noncovering parent*) of $\alpha$ if $v$ is a neighbor (parent) of at least one $a \in \alpha$ but fails to be a neighbor (parent) of at least one other $a' \in \alpha$. The set of covering neighbors (covering parents) of $\alpha$ is denoted by $\mathrm{cnb}_G(\alpha)$ [$\mathrm{cpa}_G(\alpha)$]. The set of noncovering neighbors of $\alpha$ is denoted by $\mathrm{ncnb}_G(\alpha)$. Note that $\mathrm{nb}_G(\alpha) = \mathrm{cnb}_G(\alpha) \dot\cup \mathrm{ncnb}_G(\alpha)$ and $\mathrm{cnb}_G(\mathrm{cnb}_G(\alpha)) \supseteq \alpha$. If $\alpha \equiv \{a\}$ is a singleton, then trivially $\mathrm{ncnb}_G(\alpha) = \varnothing$.

LEMMA 5.1. *Let $G \equiv (V, E)$ be a chain graph, $\xi \in \Xi(G)$ a nontrivial chain component, and $\bar{\xi} = \mathrm{cl}_G(\xi) \equiv \xi \dot\cup \mathrm{pa}_G(\xi)$. The following two properties of $G_{\bar{\xi}}$ are equivalent:*

S: *For every nonempty complete subset $\alpha \subsetneq \xi$ in $G_\xi$ such that $\kappa \equiv \mathrm{cnb}_{G_\xi}(\alpha)$ is nonempty, let $\kappa_1, \ldots, \kappa_r$ denote the connected components of $G_\kappa$. Then for each $q = 1, \ldots, r$, either:*

    (i) $\mathrm{nb}_{G_\xi}(\kappa_q) \setminus \alpha \neq \varnothing$, *or*
    (ii) $\mathrm{pa}_{G_{\bar\xi}}(\alpha) \setminus \mathrm{cpa}_{G_{\bar\xi}}(\kappa_q) \neq \varnothing$ *(equivalently, $G_{\bar\xi}$ contains a flag $t \to u - v$ with $t \in V \setminus \xi$, $u \in \alpha$, and $v \in \kappa_q$).*

S': *For every nonempty connected subset $\kappa \subsetneq \xi$ such that $\alpha \equiv \mathrm{cnb}_{G_\xi}(\kappa)$ is nonempty and complete, either:*

    (i') $\mathrm{ncnb}_{G_\xi}(\kappa) \neq \varnothing$, *or*

(ii′) $\mathrm{pa}_{G_{\bar{\xi}}}(\alpha) \setminus \mathrm{cpa}_{G_{\bar{\xi}}}(\kappa) \neq \varnothing$ *(equivalently, $G_{\bar{\xi}}$ contains a flag $t \to v - w$ with $v \in \alpha$ and $w \in \kappa$).*

PROOF. S $\Longrightarrow$ S′: Assume that $G_{\bar{\xi}}$ satisfies S and let $\kappa \subsetneq \xi$ be a nonempty connected subset in $G_\xi$ s.t. $\alpha \equiv \mathrm{cnb}_{G_\xi}(\kappa)$ is nonempty and complete. Then $\kappa' \equiv \mathrm{cnb}_{G_\xi}(\alpha) \supseteq \kappa$ is nonempty, so by S, either (i) $\mathrm{nb}_{G_\xi}(\kappa'_q) \setminus \alpha \neq \varnothing$ or (ii) $\mathrm{pa}_{G_{\bar{\xi}}}(\alpha) \setminus \mathrm{cpa}_{G_{\bar{\xi}}}(\kappa'_q) \neq \varnothing$, where $\kappa'_q$ is the unique connected component of $G_{\kappa'}$ that contains $\kappa$. If $\kappa'_q = \kappa$, then either $\mathrm{nb}_{G_\xi}(\kappa) \setminus \alpha \neq \varnothing$ or $\mathrm{pa}_{G_{\bar{\xi}}}(\alpha) \setminus \mathrm{cpa}_{G_{\bar{\xi}}}(\kappa) \neq \varnothing$, so (i′) or (ii′) holds. If $c \in \kappa'_q \supsetneq \kappa$, let $\pi = (c \equiv c_0, c_1, \ldots, c_l)$ be a minimal-length path from $c$ to $\kappa$ in $G_{\kappa'_q}$. Then $c_{l-1} \in \mathrm{nb}_{G_\xi}(\kappa) \cap \kappa'_q \subseteq \mathrm{nb}_{G_\xi}(\kappa) \setminus \alpha = \mathrm{ncnb}_{G_\xi}(\kappa)$, so (i′) holds. Thus, S′ holds in either case.

S′ $\Longrightarrow$ S: Assume that $G_{\bar{\xi}}$ satisfies S′ and let $\alpha \subsetneq \xi$ be a nonempty complete subset in $G_\xi$ s.t. $\kappa \equiv \mathrm{cnb}_{G_\xi}(\alpha)$ is nonempty. Let $\kappa_1, \ldots, \kappa_r$ denote the connected components of $G_\kappa$, fix $q \in \{1, \ldots, r\}$, and let $\alpha' \equiv \mathrm{cnb}_{G_\xi}(\kappa_q) \supseteq \alpha \neq \varnothing$. If $\mathrm{nb}_{G_\xi}(\kappa_q) \setminus \alpha \neq \varnothing$, then (i) holds, so assume that $\mathrm{nb}_{G_\xi}(\kappa_q) \setminus \alpha = \varnothing$. In this case $\alpha' = \alpha$, so $\alpha'$ is nonempty and complete, hence S′ implies that either (i′) $\mathrm{ncnb}_{G_\xi}(\kappa_q) \neq \varnothing$ or (ii′) $\mathrm{pa}_{G_{\bar{\xi}}}(\alpha') \setminus \mathrm{cpa}_{G_{\bar{\xi}}}(\kappa_q) \neq \varnothing$. If (i′), then $\mathrm{nb}_{G_\xi}(\kappa_q) \setminus \alpha \equiv \mathrm{ncnb}_{G_\xi}(\kappa_q) \neq \varnothing$ so (i) holds, while if (ii′), then (ii) is immediate, so S holds in either case. $\square$

LEMMA 5.2. *Let $\sigma$ be a strong chain component of the AMP essential graph $\mathcal{G}^*$ and let $\bar{\sigma} = \mathrm{cl}_{\mathcal{G}^*}(\sigma) \equiv \sigma \,\dot\cup\, \mathrm{pa}_{\mathcal{G}^*}(\sigma)$. Then $\mathcal{G}^*_{\bar{\sigma}}$ satisfies $\mathrm{S} \equiv \mathrm{S}'$.*

PROOF. To show that $\mathcal{G}^*_{\bar{\sigma}}$ satisfies S, let $\alpha \subsetneq \sigma$ be a complete subset of $\mathcal{G}^*_\sigma$ such that $\kappa \equiv \mathrm{cnb}_{\mathcal{G}^*_\sigma}(\alpha) \neq \varnothing$ and let $\kappa_1, \ldots, \kappa_r$ be the connected components of $\mathcal{G}^*_\kappa$. For any $q = 1, \ldots, r$, let $H$ be the graph obtained from $\mathcal{G}^*$ by replacing the lines $a - c$ by $a \to c$ for every $(a, c) \in \alpha \times \kappa_q$. This cannot create any immoralities in $H$ since $\alpha$ is complete. If $\mathrm{pa}_{\mathcal{G}^*_{\bar{\sigma}}}(\alpha) \setminus \mathrm{cpa}_{\mathcal{G}^*_{\bar{\sigma}}}(\kappa_q) = \varnothing$, this would not destroy any triplexes in $\mathcal{G}^*$. (Note that such a destroyed triplex must be a flag in $\mathcal{G}^*$, since no arrows of $\mathcal{G}^*$ are altered.) If $\mathrm{nb}_{\mathcal{G}^*_\sigma}(\kappa_q) \setminus \alpha = \varnothing$, neither would any flags or semi-directed cycles be created in $H$. Thus, $H$ would have the same skeleton and triplexes as $\mathcal{G}^*$, so $H \in \mathcal{G}$. But $H_\sigma$ contains at least one arrow, so this contradicts the assumption that $\sigma$ is strong. $\square$

We are ready to present a complete characterization of general AMP chain graphs. The following definitions are needed.

DEFINITION 5.1. Let $G \equiv (V, E)$ be a chain graph and let $\Xi_s(G)$ denote the set of nontrivial chain components $\xi$ of $G$ such that $G_{\bar{\xi}}$ contains at least one chordless undirected cycle or flag. Call a line $a - b \in G$ *strong in $G$* (indicated as $a \stackrel{\mathrm{s}}{-} b \in G$) if $a, b \in \xi$ for some $\xi \in \Xi_s(G)$. Let $\Xi_t(G)$ denote the set of trivial ($\equiv$ singleton) chain components of $G$, set $\Xi_w(G) = \Xi(G) \setminus (\Xi_s(G) \dot\cup \Xi_t(G))$, and call a line $a - b \in G$ *weak in $G$* (indicated as $a \stackrel{w}{-} b \in G$) if $a, b \in \xi$ for some $\xi \in \Xi_w(G)$. Thus, for each $\xi \in \Xi_w(G)$, $G_{\bar{\xi}}$ contains no chordless undirected cycle or flag.

DEFINITION 5.2. An arrow $a \to b \in G$ is *well protected in $G$* if it occurs in at least one of the eight configurations in Figure 6 as an induced subgraph of $G$, where $\stackrel{s}{-}$ is now defined as in Definition 5.1.

By Proposition 4.1(b), these definitions of strong/weak lines in $G$ and well-protected arrows in $G$ agree with the previous definitions when $G \equiv \mathcal{G}^*$ is an AMP essential graph.

THEOREM 5.1. *A graph $G \equiv (V, E)$ is an AMP essential graph, that is, $G = \mathcal{G}^*$ for some AMP Markov equivalence class $\mathcal{G}$, if and only if $G$ satisfies the following three conditions:*

G1: *$G$ is a chain graph, that is, is adicyclic.*
G2: *For each $\xi \in \Xi_s(G)$, $G_{\bar{\xi}}$ satisfies property $S \equiv S'$.*
G3: *Each arrow in $G$ is well protected in $G$.*

PROOF. "Only if": If $G \equiv \mathcal{G}^*$ is an AMP essential graph, G1 follows from Theorem 3.2, G2 follows from Proposition 4.1(b) and Lemma 5.2, and G3 follows from Proposition 4.3.

"If": Assume that $G$ satisfies G1, G2 and G3. Let $\mathcal{G}$ be the AMP Markov equivalence class containing $G$. To show that $G = \mathcal{G}^*$, it suffices to establish the stronger fact that:

(a) if $a \stackrel{w}{-} b \in G$, then $a \stackrel{w}{-} b \in \mathcal{G}^*$;
(b) if $a \stackrel{s}{-} b \in G$, then $a \stackrel{s}{-} b \in \mathcal{G}^*$;
(c) if $a \to b \in G$, then $a \to b \in \mathcal{G}^*$.

(a): It suffices to show that $\exists\ H, H' \in \mathcal{G}$ such that $a \to b \in H$ and $a \leftarrow b \in H'$. Since $a \stackrel{w}{-} b \in G$, $a, b \in \xi$ for some $\xi \in \Xi_w(G)$. Since $G_\xi$ is chordal, use MCS starting at $a$ (resp., $b$) to obtain a perfect orientation of $G_\xi$, thereby replacing $G$ by a graph $H$ (resp., $H'$) with $a \to b \in H$ (resp., $a \leftarrow b \in H'$). Since $G_{\bar{\xi}}$ has no flags, it is straightforward to show that $H$ and $H'$ are adicyclic and have the same triplexes as $G$, as required.

(b): It suffices to show that, for every $\xi \in \Xi_s(G)$ and every $H \in \mathcal{G}$, $H_\xi$ has no arrows. Suppose to the contrary that $H_\xi$ has at least one arrow. (Note that $H_\xi$ must be adicyclic since $H$ is adicyclic.) Choose $a \to b \in H_\xi$ so that $b$ is maximal w.r.t. the pre-ordering induced on $\xi$ by $H_\xi$, that is, there exists no semi-directed path in $H_\xi$ beginning at $b$. Let $\xi(b) \in \Xi(H_\xi)$ be the unique chain component of $H_\xi$ (possibly trivial) containing $b$ and let $\alpha = \mathrm{nb}_{G_\xi}(\xi(b))$. Note that $a \in \alpha$ since $a - b \in G_\xi$. For any $a' \in \alpha$, $a' - b' \in G_\xi$ for some $b' \in \xi(b)$, so $a' \cdots b' \in H_\xi$. If $a' \neq a$, then, by the maximality of $b$ and the connectedness of $\xi(b)$, $a' \Rightarrow b' \in H_\xi$, so by the definition of $\xi(b)$, $a' \to b' \in H_\xi$. If $a' = a$, then $a' \to b \in H_\xi$. Thus, for every $a' \in \alpha$, $a' \to b'' \in H_\xi$ for all $b'' \in \xi(b)$ since $\xi(b)$ is connected and $H_\xi$ can have no flags (because $G_\xi$ has no arrows, hence no triplexes). Thus, $\alpha = \mathrm{cnb}_{G_\xi}(\xi(b))$, which in turn implies that $\alpha$ is complete, since $H_\xi$ can have no immoralities. Therefore, since $G_{\bar{\xi}}$ satisfies property S' by G2, either $\mathrm{ncnb}_{G_\xi}(\xi(b)) \neq \varnothing$ or $\mathrm{pa}_G(\alpha) \setminus \mathrm{cpa}_G(\xi(b)) \neq \varnothing$. The former is impossible since $\mathrm{ncnb}_{G_\xi}(\xi(b)) = \mathrm{nb}_{G_\xi}(\xi(b)) \setminus \mathrm{cnb}_{G_\xi}(\xi(b)) = \alpha \setminus \alpha = \varnothing$. If the latter holds, then $\exists v \in V \setminus \xi$, $a' \in \alpha$, and $b' \in \xi(b)$ such that $v \to a' \in G$ but $v \cdot / \cdot b'$ in $G$. But necessarily $a' - b' \in G_\xi$, while $a' \to b' \in H_\xi$, so $v \to a' - b'$ occurs as a flag in $G$ but not in $H$, also a contradiction.

(c): It suffices to show that
$$A \equiv \{a' \in V \mid \exists b' \in V \ni a' \to b' \in G, a' - b' \in \mathcal{G}^*\} = \varnothing.$$
If $A \neq \varnothing$, let $a$ be a minimal element of $A$ with respect to the pre-ordering induced on $V$ by $G$. Since $a \in A$,
$$B \equiv \{b' \in V \mid a \to b' \in G, a - b' \in \mathcal{G}^*\} \neq \varnothing.$$
Let $b$ be a minimal element of $B$; in particular, $a \to b \in G$ and $a - b \in \mathcal{G}^*$. By G3, $a \to b$ is well protected in $G$, so it occurs in at least one of the eight configurations (i)–(viii) in Figure 6 as an induced subgraph of $G$.

(i) If $c \to a \to b$ occurs as an induced subgraph in $G$, then by the minimality of $a$ the flag $c \to a - b$ must occur in $\mathcal{G}^*$, contradicting the fact that $G$ and $\mathcal{G}^*$ have the same triplexes.

(ii) If $a \to b \leftarrow c$ occurs as an induced subgraph in $G$, then the flag $a - b \leftarrow c$ must occur in $\mathcal{G}^*$, which would require that $a \stackrel{\mathrm{s}}{-} b \in \mathcal{G}^*$, hence $a - b \in G$, again a contradiction.

(iii) If $a \to b \leftarrow c \leftarrow a$ occurs as a triangle in $G$, then, by the minimality of $b$, the semi-directed triangle $a - b \Leftarrow c \leftarrow a$ must occur in $\mathcal{G}^*$, contradicting the adicyclicity of $\mathcal{G}^*$.

(iv) If $a \to b$ occurs in configuration (iv) in $G$, denote the immorality in this configuration by $c \to b \leftarrow d$. Since $G$ and $\mathcal{G}^*$ have the same triplexes, either $c \to b \Leftarrow d$ or $c \Rightarrow b \leftarrow d$ occurs as an induced subgraph of $\mathcal{G}^*$. Without loss of generality, assume the former. Since $a - b \in \mathcal{G}^*$, necessarily $c \to a \in \mathcal{G}^*$ by the adicyclicity of $\mathcal{G}^*$, hence $a \to d \in \mathcal{G}^*$ since $G$ and $\mathcal{G}^*$ have the same triplexes. But then the semi-directed triangle $a \to d \Rightarrow b - a$ occurs in $\mathcal{G}^*$, again contradicting adicyclicity.

(v), (vi), (vii), (viii): If $a \to b$ occurs in configuration (v), (vi), (vii) or (viii) in $G$, by (b) the strong line in this configuration must also be strong in $\mathcal{G}^*$. Since $a - b \in \mathcal{G}^*$, Lemma 4.1(d) implies that $a \stackrel{\mathrm{s}}{-} b \in \mathcal{G}^*$, which contradicts $a \to b \in G$. This completes the proof of (c). □

The necessity of Condition G2 is demonstrated by the graph $G$ in Figure 2 and the upper graph $G_\infty$ in Figure 1. The necessity of G3 is demonstrated by the graph consisting of a single arrow.

For comparison with Theorem 5.1, the characterization of ADG essential graphs in Theorem 4.1 of [3] can be stated as follows: $G$ is an ADG essential graph iff it satisfies conditions G1, G2': $\Xi_s(G) = \varnothing$ and G3.

More detailed results on the structure of a strong chain component $\sigma$ of $\mathcal{G}^* \equiv (V, E^*)$ can be obtained from Theorem 5.1. For example, if $\mathcal{G}^*$ is a strong connected undirected AMP essential graph (so $\sigma = V$), then by Proposition 4.1(b), $\mathcal{G}^*$ must contain at least one chordless cycle. Using Theorem 5.1, however, it can be shown that unless $\mathcal{G}^*$ consists exactly of a single chordless cycle, it must contain at least one additional chordless cycle. For a general AMP essential graph $\mathcal{G}^*$, by Proposition 4.1(b), $\mathcal{G}^*_{\bar\sigma}$ must contain either a chordless undirected cycle or a biflag. Using Theorem 5.1, however, it can be shown that $\mathcal{G}^*_{\bar\sigma}$ must contain another chordless undirected cycle or biflag or "3-halfbiflag." These results appear in [5, 6].

**6. Current research.** Like the essential graph $D^*$ for ADG Markov models, the AMP essential graph $\mathcal{G}^*$ plays a fundamental role for inference, model selection and model averaging for AMP CG Markov models. For these purposes, the results of [3] and [21] can be extended to AMP CG models by means of our characterization of AMP essential graphs in Theorem 5.1 above. In particular, a polynomial-time algorithm for constructing $\mathcal{G}^*$ from any $G \in \mathcal{G}$ has been obtained—see [5], Chapter 7.

The following additional topics are currently under development:

1. The computational complexity of the above construction algorithm.
2. An algorithm for recovering all $G \in \mathcal{G}$ from $\mathcal{G}^*$.
3. Markov chain Monte Carlo algorithms for model search over the space of AMP Markov equivalence classes by means of AMP essential graphs.
4. A catalog of AMP essential graphs with small vertex sets.
5. Determination of the ratio $r_n$ of the number of AMP CG Markov equivalence classes to the total number of CGs with $n$ vertices (cf. [16] for the corresponding question for ADG models).
6. When is an AMP CG Markov equivalent to some LWF CG, and vice versa? (Cf. [4], Theorem 6.)

**Acknowledgments.** We thank Robert Castelo, Sanjay Chaudhuri, Mathias Drton, Steven Gillispie, David Madigan, Chris Meek, Marina Meila, Thomas Richardson, Alberto Roverato, Milan Studený, a referee and especially an Associate Editor for helpful comments and suggestions.

DEPARTMENT OF MATHEMATICS
INDIANA UNIVERSITY
BLOOMINGTON, INDIANA 47405
USA
E-MAIL: standers@indiana.edu

DEPARTMENT OF STATISTICS
UNIVERSITY OF WASHINGTON
SEATTLE, WASHINGTON 98195-4322
USA
E-MAIL: michael@stat.washington.edu